# A non-damped stabilization algorithm for multibody dynamics

Igor Fernández de Bustos · Haritz Uriarte · Gorka Urkullu · Vanessa García-Marina

**Abstract** The stability of integrators dealing with high order Differential Algebraic Equations (DAEs) is a major issue. The usual procedures give rise to instabilities that are not predicted by the usual linear analysis, rendering the common checks (developed for ODEs) unusable. The appearance of these difficult-to-explain and unexpected problems leads to methods that arise heavy numerical damping for avoiding them. This has the undesired consequences of lack of convergence of the methods, along with a need of smaller stepsizes. In this paper a new approach is presented. The algorithm presented here allows us to avoid the interference of the constraints in the integration, thus allowing the linear criteria to be applied. In order to do so, the integrator is applied to a set of instantaneous minimal coordinates that are obtained through the application of the null space. The new approach can be utilized along with any integration method. Some experiments using the Newmark method have been carried out, which validate the methodology and also show that the method behaves in a predictable way if one considers linear stability criteria.

I. Fernández de Bustos · H. Uriarte (✉) · G. Urkullu
Department of Mechanical Engineering, Faculty of Engineering of Bilbao, University of the Basque Country, Alameda Urquijo s/n, 48013 Bilbao, Spain
e-mail: haritz.uriarte@ehu.eus

I. Fernández de Bustos
e-mail: igor.fernandezdebustos@ehu.eus

G. Urkullu
e-mail: gorka.urkullu@ehu.eus

V. García-Marina
Department of Mechanical Engineering, University School of Engineering of Vitoria-Gasteiz, University of the Basque Country, Nieves Cano 12, 01006 Vitoria-Gasteiz, Spain
e-mail: vanessa.garcia@ehu.eus

## 1 Introduction

### 1.1 Background

Ordinary Differential Equation systems (ODEs) represent a great number of physical phenomena including structural analysis, multibody dynamics, thermal conduction or electrical systems. This has led to intensive investigation about them [1–6]. One of the key issues in ODEs is the stability. If the system to integrate is stiff, the integrator might fail due to instability. A-Stable integrator are stable regardless of the stepsize, which is quite interesting. Explicit methods, which can deliver great performance in terms of convergence, cannot be A-Stable. Also, quite good performing implicit methods are not A-Stable.

This leads to the concept of A($\alpha$) stability, which introduces a conditional stability. This means that there is a limit in the timestep that can be used in the integration.

Often, these problems, depending on the modeling of the system, require the incorporation of some constraints. Therefore, a nonlinear Differential Algebraic Equations system (DAEs) must be solved to predict the behavior of many of these physical phenomena [7–10]. Commonly, tools of optimization have been taken into account to solve the nonlinear DAEs as it can be seen in Nocedal [11]. In any case, it seems clear that the high index systems are more difficult to solve [12] and are associated with a singular mathematical set up [13, 14]. In this way, one of the main strategies to face this type of problems is the index reduction technique which applies a mathematical process to reduce the index of a set of DAEs. A very interesting work in this field was done by Natsiavas in [15]. In this work, the authors present a new theoretical approach for deriving an appropriate set of equations of motion for a class of mechanical systems subjected to motion constraints. Another fascinating approach was presented by González [16] and Bayo [17], who used the penalty factors present in the augmented Lagrange formulation to determine the reaction forces in redundant constrained dynamic systems in order to represent the physical properties of the system in the model. However, index reduction through the analytical differentiation of the constraint equations causes the progressive drift of the computed solution; therefore, it is necessary to apply some kind of stabilization. To alleviate this problem Bayo and Cuadrado [18, 19] propose a post-stabilization technique based on coordinates projection. In addition, apart from the index reduction, there is another well-known approach: the coordinates reduction technique. These methods try to divide the coordinates in two groups: dependent and independent coordinates. An excellent work was done by Zhang in [20] using this technique. He uses an implicit Runge–Kutta method for the solution of index-3 DAEs. The iteration of the constraint equations is embedded in the iteration of the non-linear algebraic equations that originate from the implicit method and, using the coordinate partitioning technique, the independent coordinates of the set of coordinates of the system are chosen. An interesting study of the behavior of these two main techniques can be found in the work done by Jalón [21]. He contributed to the resolution and clarification of the multiplicity of solutions to the constraint equations, where they focused on three methods: the Lagrange equations of the first kind, the null space method (introduced by Schwerin [22]), and the Maggi equations. In this sense, Arnold et al. [23] studied the effect of velocity projections on stability, and Cuadrado et al. [24] performed a comparison between four methods to simulate multibody dynamics with constraints. These methods were the augmented Lagrange formulation (ALF) index-1 and index-3 with projections, ALF-1 and ALF-3, respectively, a modified state-space formulation, and a fully recursive formulation.

### 1.2 Formulation of the problem of interest for this investigation

In the case of ODEs, stability criteria are quite well established. But, unfortunately, these criteria cannot be directly translated to DAEs when applying ODE integrators in the usual way. Gear in [25], highlighted the difficulties associated with the solution of this type of systems. Higueras et al. in [26] analyze the stability of index-2 DAEs, remarking the unexpected stepsize limitations that appear even when an A-stable method is used. In this paper, the authors conclude that choosing a suitable formulation of DAEs is a vital task for a successful and effective numerical integration as much as the numerical method. Otherwise, stepsize limitations might appear because implicit methods behave in the same way as the explicit ones. Even worse, the method may fail completely. These problems were already observed by Ascher in [27], where the implicit Euler method was transformed into an explicit method when it was applied to solve some problems. About stability in linear index-2 DAEs, Hanke [28] presents a study of asymptotic properties of solutions on infinite intervals. The author states, as it happens in index-1 DAEs [29], that algebraically stable implicit integration methods are shown to be B-stable, provided that the null space of the leading Jacobian is constant. If this null space rotates, stability properties may change. Many others have also taken note in these problems, as it can be seen in the work done by Liu [30] which studies the stability of the numerical methods for linear index-3 DAEs.

As it can be seen, the stability of the integrator methods is a significant problem, especially in constrained high order differential equation systems. The assumptions adopted in linear problems are not applicable to nonlinear ones and, usually, instabilities appear. These undesirable problems lead to the application of an integration method that introduces a certain numerical damping (L-Stable methods). Nevertheless, the convergence towards the solution of the system is broadly affected and, obviously, it forces to use smaller time steps.

Another quite remarkable family of solutions to this problem is the use of a reduced set of minimal coordinates (see [31–33]). In these methods, the authors select a set of the original coordinates and the integration is performed in this set. Afterwards, the restrictions are used to obtain the values of the rest of the coordinates. Unfortunately, this reduced set cannot usually be applied to the whole integration interval and, thus, the set has to be changed from time to time. This leads to an increase in the computational cost, which has commonly discouraged the use of these methods.

An interesting aspect to point out is that explicit methods, while being $A(\alpha)$-stable, do not show unpredictable behavior. This means that the timestep required for stability can easily be predicted with classical ODE theory. For example, methods as those presented in [34] exhibit predictable behavior. The main problem of these methods is that their stability limits the timestep required to solve the problem. In the case of some stiff problems, this leads to the need of a high computational cost derived just from the stability requirement.

### 1.3 Scope and contribution of this study

In this paper a new approach for applying numerical integration to index 3 DAEs (although easily portable to DAEs of other indexes) is presented. The main idea of the algorithm is to integrate in the tangent space to the manifold. The new approach presents a predictable stability behavior and can be applied along with any kind of integration method such as Newmark, HHT and others. The method also allows one to force the simultaneous verification of constraints in terms of the function itself and all the relevant derivatives. Although similar in concept to the methods based in a reduced set of coordinates, this method is quite different in the sense that the coordinates are different in each step and represent a movement in the tangent space. In order to do this without hampering computational cost, most of the needed information is obtained from the algorithm employed to solve the equations. Some numerical examples are presented which show the behavior of the algorithm in stiff problems. Also, a quite simple explanation of the phenomenon of unpredictable stability in DAEs in the so-called structural integrators is presented.

### 1.4 Organization of the paper

The paper is organized as follows. First, the usual approach for applying integration schemes to index 3 DAEs is presented. Second, a quite graphical example of the problems that might arise in stiff problems due to stability issues is introduced. Afterwards, an explanation of the usual problems that appear in the commonly used approach employed for index 3 DAEs integrators is given. The next point addresses the new approach for applying the integrators. Some numerical experiments are presented to demonstrate the new approach behavior and, finally, some conclusions are drawn.

## 2 Common approaches for the integration of constrained systems with structural integrators

We consider high index DAEs in the form expressed by the following equations:

$$F(\dot{x}, x, t) = 0 \quad (1)$$

$$q(t, x) = 0 \quad (2)$$

where (1) is the differential equation to be solved, and (2) are the boundary conditions. In the case of index-3 DAEs, which are typical of Multibody Dynamics, Structural Nonlinear Mechanics and other engineering problems, one can write Eq. (1) as:

$$\boldsymbol{M}(\boldsymbol{x})\ddot{\boldsymbol{x}} = \boldsymbol{f}(\boldsymbol{x}, \dot{\boldsymbol{x}}, t) + \boldsymbol{G}^T(\boldsymbol{x})\boldsymbol{\lambda} \quad (3)$$

where $\boldsymbol{M}(\boldsymbol{x})$ is the mass matrix, which is invertible in the integration interval, $\boldsymbol{f}(\boldsymbol{x}, \dot{\boldsymbol{x}}, t)$ are the applied, Coriolis and velocity-dependent forces, and $\boldsymbol{G}^T(\boldsymbol{x})\boldsymbol{\lambda}$ are the constraint forces. Two main approaches are

374    Meccanica (2022) 57:371–399

considered to apply an ODE integrator to this problem. The traditional approach reformulates the problem, resulting in a system in the form expressed in the following equations and Eq. (2).

$$\dot{x} = u \tag{4}$$

$$M(x)\dot{u} = f(x, u, t) + G^T(x)\lambda \tag{5}$$

With this procedure, one can apply a first order method, which can be implicit (usually single step) or explicit (usually multi step).

A more recent approach takes advantage of the so called structural integrators ([34–37]). These methods are directly applied to the system expressed by eqns. (3) and (2). The main advantage of these methods is that they do not duplicate the number of functions to integrate. In any case, one can also find implicit (for example: Newmark, HHT) and explicit (for example Central Differences) methods. Again, usually implicit methods are formulated in a single step approach and explicit methods are usually multi stepped.

In any case, usually one applies the integrator to Eq. (3) (in the structural integrator approach) or to eqns. (4) and (5) (in the first order approach) and, afterwards, the constraints presented in Eq. (2) are introduced. We will further develop the algorithm for the structural approach, although a similar development can be done for the first order approach.

The first step is to obtain a linearization of Eq. (3). For example, with an implicit method, the linearization would be usually performed in $(t + \Delta t)$, and thus one would reach:

$$M_L\ddot{x}(t + \Delta t) + C_L\dot{x}(t + \Delta t) + K_Lx(t + \Delta t) \\ = f_L + G_L^T\lambda \tag{6}$$

Usually the terms $M_L$, $C_L$, $K_L$, $G_L$ and $f_L$ are obtained in a numerical way. These terms differ from their counterparts $M$, $C$, $K$, $G$ and $f$ and are the result of applying a linearization such as Taylor series to an otherwise non-linear Eq. (3). Terms $M_L$, $C_L$, $K_L$, $G_L$ and $f_L$ therefore relate not only to $M$, $C$, $K$, $G$ and $f$ respectively; as a result of the applied linearization (and their dependencies on the system's coordinates) different derivatives of the same term with respect to different coordinates arise, forming several terms that are each a function of $x$, $\dot{x}$, $\ddot{x}$ or $\lambda$. All these dependencies related to each one of the system's variables are summed together separately; this means that, for example for $x$, all the dependencies related to $x$ form $K_L$ as the matrix that includes all the linearized dependencies related to $x$ and is multiplied by the corresponding system's variable in the equilibrium equation, as seen in (6). The same happens with $C_L$ related to velocities, to $f_L$ related to independent terms, and to the remaining variables obtained as a result of the linearization that form Eq. (6). The equation obtained with this linearization process requires an iterative scheme to obtain the correct solution.

One can use the function to be integrated to solve the problem or reduce the index. Thus, if the integration is based on the function itself, the procedure is referred as Index 3 Formulation. In this case, after applying the integrator equations to Eq. (1), one reaches an equation (expressed for implicit methods) in the form of:

$$A_1 x(t + \Delta t) = g_1 + G_{L1}^T \lambda \tag{7}$$

If, instead, one uses the first derivative (Index 2 Formulation), one reaches:

$$A_2 \dot{x}(t + \Delta t) = g_2 + G_{L2}^T \lambda \tag{8}$$

Finally, using the second derivative (Index 1 Formulation), one would reach:

$$A_3 \ddot{x}(t + \Delta t) = g_3 + G_{L3}^T \lambda \tag{9}$$

Being $A_1, A_2, A_3, G_{L1}, G_{L2}, G_{L3}, g_1, g_2, g_3$ matrices and vectors that depend on the current estimation of $x(t + \Delta t), \dot{x}(t + \Delta t)$ or $\ddot{x}(t + \Delta t)$.

The linearization of Eq. (2) will lead to a system of equations that would allow one to solve eqs. (7), (8), or (9). Obviously, the linearization would be different if, instead of the function, one of its derivatives is employed. For example, using an Index 3 formulation, one would obtain from Eq. (2) an expression in the following form:

$$H_L x(t + \Delta t) = b_L \tag{10}$$

where $H_L$ is the matrix that includes the linearized expressions of the applied constraints that are multiplied by the solution vector expressed in the system's coordinates, and $b_L$ is the independent term vector obtained after the same linearization that provided $H_L$. This linearization will be explained in Sect. 4. Similarly, matrix $G_L$ used previously in this document is the matrix that includes the linearized expressions of the applied constraints expressed in the coordinates in which the forces are defined, which can differ from the ones employed to form $H_L$ as Urkullu explains [34]. For instance, $G_L$ requires rotations for the angular coordinates, and $H_L$ can be formed using any other angular coordinates that can be considered more convenient, such as quaternions. Both matrices have been named differently to avoid confusion related to this matter despite the fact that both define the Jacobian of the system.

Thus, the system to be solved is composed of eqs. (7) and (10), which are modified in an iterative approach. With this in mind, the algorithm is quite simple. Figure 1 presents the mentioned algorithm when formulated in terms of the function expressed before.

If one uses Index 2 or Index 1 formulations, a drift-off may appear and this leads to the need of the so called stabilization methods, such as Baumgarte or Projection. ([38]).

One could consider a DAE as differential equations on a manifold (see [39]). In this sense, this way of integrating the DAE could be explained as follows: integrating the function while forcing the solution to belong to the manifold.

There are some criticisms that can be applied to this approach. For instance, one of them is the fact that one cannot impose that the function to be integrated and its derivatives satisfy the constraints in the considered integration step. This is easy to demonstrate. Let us consider the particular case of a system with N-functions to integrate, N differential equations and $N_R$ constraints, being integrated using Newmark. This is a common case in multibody dynamics. The constraints bring about $N_R$ Lagrange multipliers. Thus, for each integration step there are $3N + N_R$ unknowns to solve.

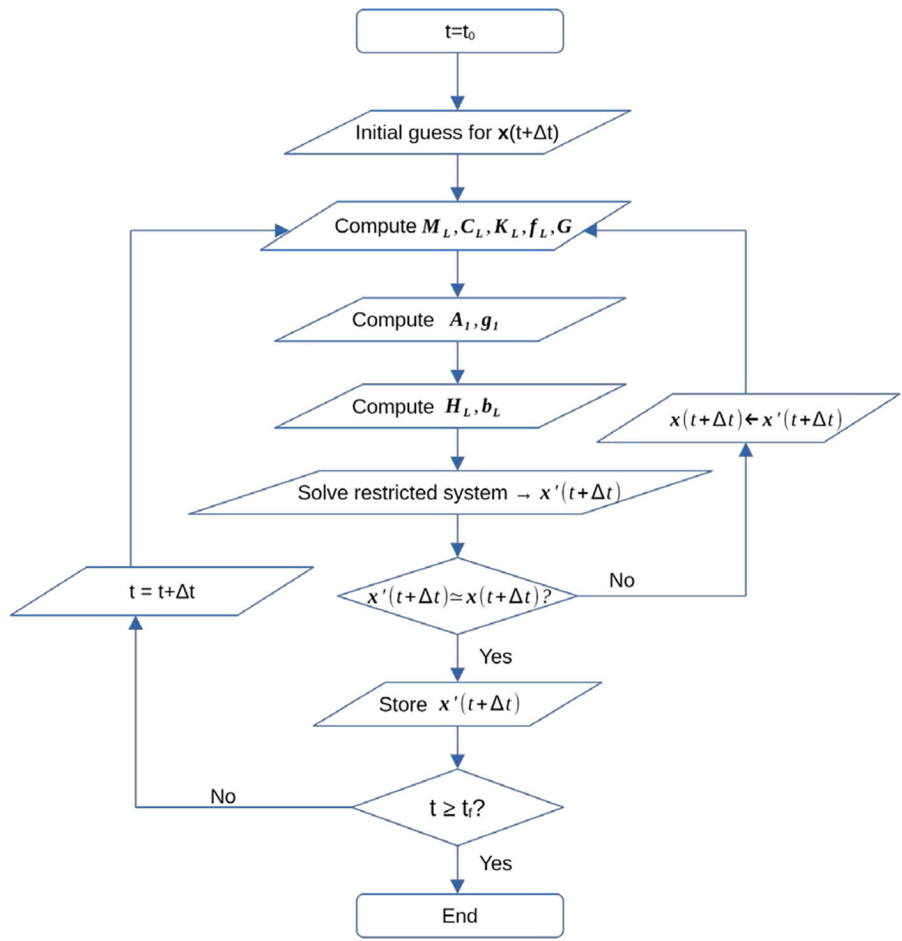

**Fig. 1** Typical integration algorithm

For this, one has 2 N Newmark equations, the equilibrium in $t + \Delta t$ leads to another set of N equations and, finally, one gets $3N_R$ equations from the constraints formulated for the functions and their two derivatives. Thus, one has $3N + 3N_R$ equations for $3N + N_R$ unknowns. Therefore, one must take the choice of applying constraints solely to the function ignoring the derivatives, applying the constraints solely to one of the derivatives ignoring the function and other derivatives (leading to the already commented drift problem) or applying a compromise solution. Another issue is the erratic behavior regarding stability.

The stability of ODEs is a quite thoroughly studied subject for first order methods. A-Stability is of most importance in stiff problems, because it will assure a solution (although not necessarily accurate). Explicit methods can be quite efficient, but cannot be A-Stable. In these cases one speaks of $A(\alpha)$-Stability, which means that the method is stable depending on some conditions. In the case of structural integrators, usually one speaks about conditional and unconditional stability. Conditional stability is similar to $A(\alpha)$—Stability, while unconditional stability is equivalent to A-Stability.

Some methods exhibit a configurable stability behavior. An example is the Newmark method, which depends on two parameters which are selectable, $\alpha$ and $\beta$. The Newmark method requires the following equation to be verified for an integration to be stable [40]:

$$\Delta t \leq \frac{1}{\omega_{max}} \sqrt{\frac{1}{\frac{1}{2}\alpha - \beta}} \qquad (11)$$

Being $\omega_{max}$ the highest natural frequency of the system to integrate. If $\alpha > 0.5$ and $\beta > 0.5\alpha$, the method is unconditionally stable (A-Stable), while if $\alpha > 0.5$ and $\beta < 0.5\alpha$, the method is conditionally stable. But unfortunately, when one applies ODE integrator to DAEs, the ODE stability analysis cannot be used. This is quite an issue, because it can considerably limit the step sizes that can be used and, consequently, the efficiency. Furthermore, behavior can be quite unpredictable. A common solution to this problem is the use of L-Stable methods. These methods introduce heavy damping in the high frequencies, thus stabilizing the problem. An example is the HHT method, which is currently one of the most used ([36, 37, 41]). But they also introduce damping in the lower frequencies, which sometimes can be an issue.

Another approach that can be used is to integrate the restricted system in a set of minimal coordinates obtained using Null-Space algorithms or other approaches [32, 33, 42]. These methods use the same set of coordinates along part of the integration time and, whenever is needed, they are changed. These methods usually are slow and lead to lack of convergence.

## 3 A simple example of the problem. a stiff pendulum

In order to expose the problems of the usual approach, a simple example has been designed. It consists of a simple pendulum affected by gravity (see Fig. 2). A harmonic torque with a quite low frequency is applied to it as shown in Eq. (12). This should allow one to use a considerably large stepsize, but depending on the integration method employed it can lead to instability.

$$T = T_0 sin(\overline{\omega} t) \qquad (12)$$

The system of reference is on the fixed node and gravity is considered as $g = 9.8 m/s^2$. The concentrated mass of the pendulum is $M = 1 kg$. The length of the truss is $L = 1m$. The truss weight is considered to be 0. The parameters defining the torque are $T_0 = 0.1 Nm$ and $\overline{\omega} = 0.1 rad/s$. The initial conditions are $\theta(0) = 0$, $\dot{\theta}(0) = 0$. We must point out that, in this problem, the mass matrix should be underdetermined, but this issue can be easily solved applying constraints to the system, or using minimal least squares solutions.

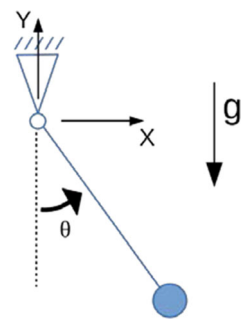

**Fig. 2** A simple pendulum

If one expresses the problem using $\theta$ as a parameter, the equilibrium equation can be expressed as shown in the following eqn.

$$mL^2\ddot{\theta}(t) + mgLsin(\theta(t)) - T(t) = 0 \quad (13)$$

A linearization of the system of equations can be obtained in an analytical way:

$$mL^2\ddot{\theta}(t) + mgLsin(\theta_0) + mgLcos(\theta_0)(\theta(t) - \theta_0) - T(t) = 0 \quad (14)$$

$$mL^2\ddot{\theta}(t) + mgLcos(\theta_0)\theta(t) = T(t) - mgLsin(\theta_0) + mgLcos(\theta_0)\theta_0 \quad (15)$$

Thus, for this system, one has $M_L = mL^2$ and $K_L = mgLcos(\theta_0)$. This means that the natural frequency of this system can be obtained from:

$$\omega = \sqrt{\frac{mgLcos(\theta_0)}{L^2}} = \sqrt{\frac{gcos(\theta_0)}{L}} \quad (16)$$

The maximum value happens in 0:

$$\omega_{max} = \sqrt{\frac{g}{L}} \quad (17)$$

which, for the presented configuration, yields: $\omega_{max} = 3.1304951 rad/s$. This is the parameter that would define the stability limits if one accepts the values obtained for linear systems. For example, for a central difference method [34], one would reach:

$$\Delta t_{maxs} = \frac{2}{\omega_{max}} = 0.6388765s \quad (18)$$

Let us consider now the step size required to correctly represent the problem, regardless of stability. If the effect of the torque is considered, providing for the need of about ten points to represent a cycle of the torque, one would use:

$$\Delta t_{maxt} = \frac{2 \cdot \pi}{10\overline{\omega}} \simeq 6.28s \quad (19)$$

Furthermore, if one wants to represent correctly the oscillations produced by gravity, one would require (in a conservative estimation) the value obtained in:

$$\Delta t_{maxg} = \frac{2 \cdot \pi}{10\omega_{max}} \simeq 0.2s \quad (20)$$

The dynamic effects of gravity are of such low amplitude that they can be disregarded. Thus, a correct step size would be about 6. But this would violate (at least in the case of the central differences method) the constraint imposed by stability criteria. Therefore, this is a good yet simple benchmark to study the stability issues. We aim to verify whether the usual way of applying constraints is to blame or not for the instability of the problem, for which we will first solve the problem without the need of employing constraint equations. In order to do so, an unrestricted minimal set of coordinates will be used. Note that this is a particular case in the sense that one can integrate the whole problem in an unrestricted system of coordinates. This is not so common. Now this equation shall be integrated using different schemes. A central difference approach and two different configurations of Newmark: the Fox-Goodwin approach and the Trapezoidal rule.

Applying a central difference scheme, one reaches the following equations for the general case and for the first instant:

$$\theta(t + \Delta t) = 2\theta(t) - \theta(t - \Delta t) + \frac{\Delta t^2}{mL^2}(T(t) - mgLsin\theta(t)) \quad (21)$$

$$\theta(0 + \Delta t) = \theta(0) + \Delta t\dot{\theta}(0) + \frac{\Delta t^2}{2mL^2}(T(0) - mgLsin\theta(0)) \quad (22)$$

In Fig. 3 the obtained results are presented for step sizes of 0.01 s, 0.1 s, 0.6 s and 0.7 s.

A brief explanation regarding nomenclature can be added here. Due to linearization the term $\theta_0$ found in Eqs. (14), (15) and beyond is the value the variable $\theta$ has in the time step that is being processed, but in the previous iteration. Consequently, $\theta_0$ is updated with every iteration. This term is not to be confused with the value this same variable has on the initial conditions, which is noted as $\theta(0)$.

The results in minimal coordinates behave as predicted by linear estimations. With step sizes lower than the theoretical stability limit, $\Delta t_{maxs}$, the method is stable and delivers reasonable results. It is also interesting to observe how the step size of 0.6 fails to deliver a good representation of the small oscillation (a step size of 0.2 was estimated as a good choice). As predicted by Eq. (18), a step size of 0.7 leads to instability.

Let us now consider the solution keeping a minimal coordinates approach, but now using a Newmark approach with a conditionally stable configuration. A 4th order Fox-Goodwin approach has been used ($\alpha = 1/2, \beta = 1/12$). This is a conditionally stable approach, and the linear stability condition can be expressed as Eq. (10). Thus, in this case: $\Delta t_{maxs} \leq 0.78246 s$. The results with stepsizes of 0.01 s, 0.1 s, 0.6 s, 0.7 s, 0.78 s and 0.79 s are presented in Fig. 4.

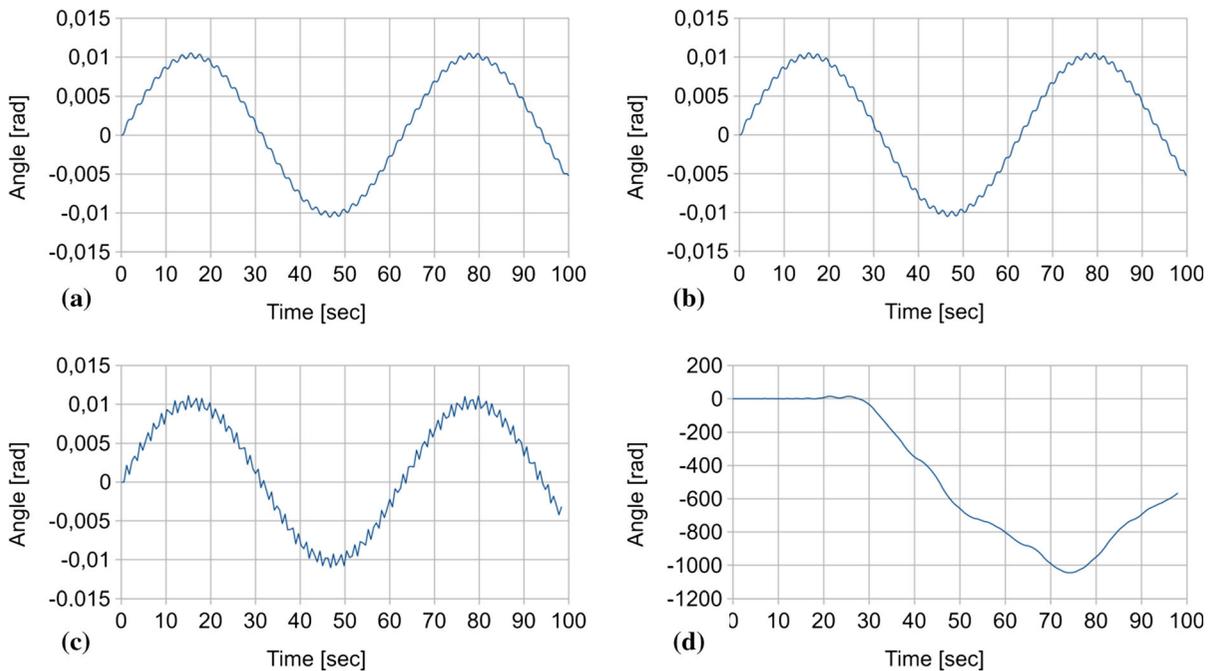

**Fig. 3** Central difference, minimal coordinates solution for step sizes **a** 0.01 s, **b** 0.1 s, **c** 0.6 s and **d** 0.7 s

Again, results are quite predictable. Step sizes of 0.7 s and below lead to good results. However, a stepsize of 0.78 s represents correctly the effect of the force, but it is too large to correctly solve a force in the frequency of the movement derived from the gravity effects. This is the reason behind the considerable oscillation of the solution. Nevertheless, it is important to point out that the system is stable. Finally, the stepsize of 0.79 s leads to unstable integration, as predicted by Eq. (11).

Let us consider now an unconditionally stable Newmark scheme, by taking $\alpha = 0.5$ and $\beta = 0.25$. This approach is commonly known as trapezoidal rule. The chosen stepsizes are: 0.1 s, 0.7 s, 0.79 s and 6 s. Results are presented in Fig. 5.

Again, results are as predicted by the linear theory. The unconditional stability of trapezoidal rule allows this scheme to report reasonable results even with a stepsize of 6 s.

Now we consider the use of a constrained system approach, for which the differential equation that provides the equilibrium equation is:

$$N_G^T M_L \ddot{x}(t + \Delta t) + N_G^T C_L \dot{x}(t + \Delta t) + N_G^T K_L x(t + \Delta t) = N_G^T f_L \tag{23}$$

This equation is the result of premultiplying (6) by the null space of $G_L$; then the constraints are introduced. The term "constrained system" is used to name those systems where more variables than the bare minimum are used to define the state of the problem (in the case of multibody dynamics the position of the elements) and thus, one needs to introduce constraints in the system of equations to solve the problem. This happens not only in all general approaches but also in some particular situations where one tries to use minimal coordinates, but part of the problem is difficult to introduce without constraints. One could name these approaches as "mixed". In this case, general coordinates $(x, y, \theta)$ will be used to define the situation of the center of gravity of the pendulum. The location of the revolute joint is fixed (2 constraints). The total number of

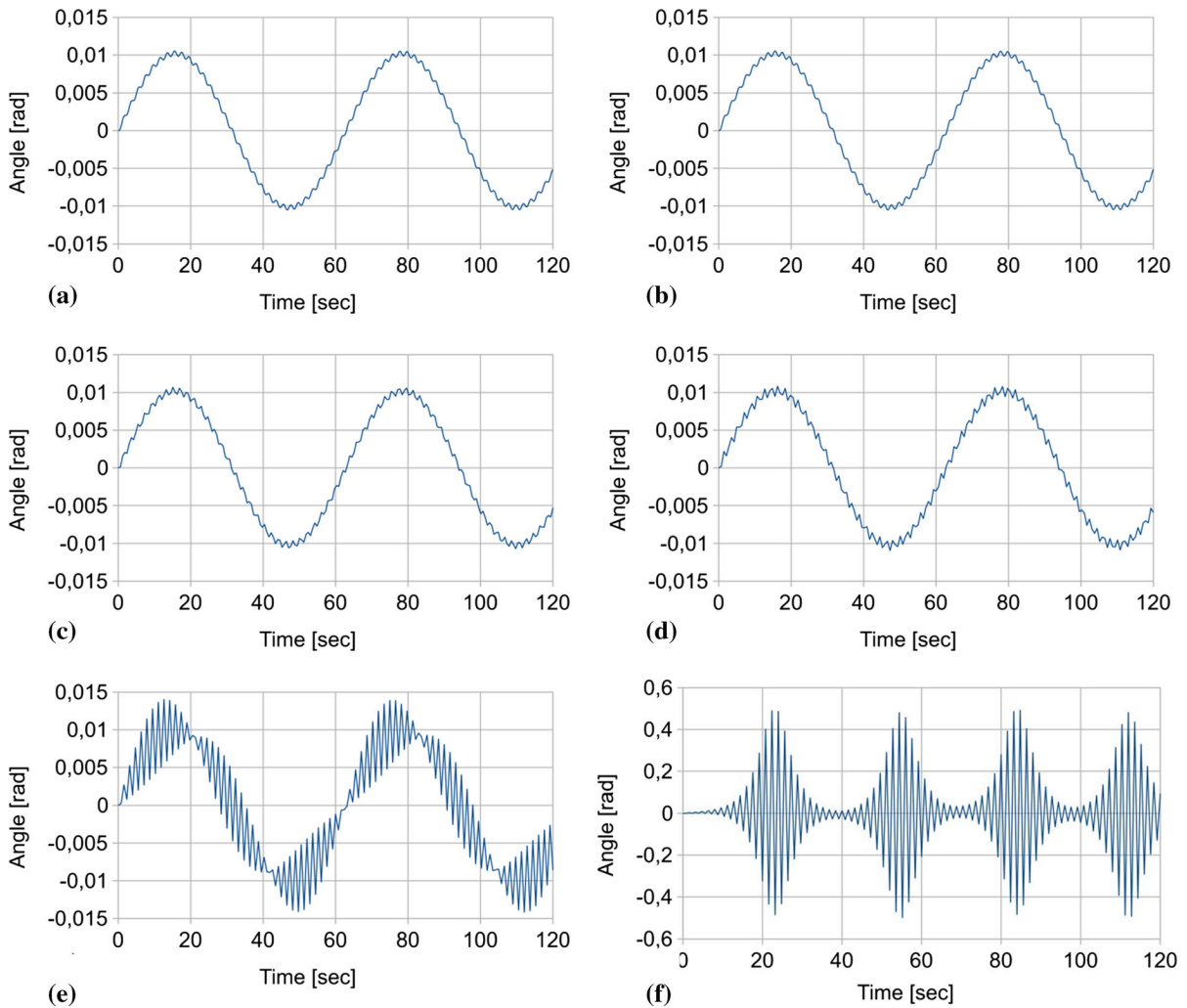

**Fig. 4** Fox-Goodwin, minimal coordinates solution for step sizes **a** 0.01 s, **b** 0.1 s, **c** 0.6 s, **d** 0.7 s, **e** 0.78 s and **f** 0.79 s

variables is 3 and the number of constraints is 2, as shown in the following equation:

$$q_R = r_G + R \cdot \bar{r}^L_{GR} - r_C = 0 \qquad (24)$$

where $q_R$ is the expression of the constraint, $r_G$ is the position of the constrained solid's center of gravity, $R$ is its rotation matrix, $\bar{r}^L_{GR}$ is the distance from the solid's center of gravity to the point where the constraint is applied expressed in local coordinates, and $r_C$ is the position of the rotation constraint in a fixed position, expressed in global coordinates. In order to solve the problem we have used a structural integrator approach, similar to those presented by Cuadrado et al. [24]. In the trapezoidal rule configuration, the method seems to work in a correct way, as can be seen in the results shown in Fig. 6.

The method seems to work well using a trapezoidal rule, but this is not as usual as one might think. There are plenty of reports of trapezoidal rule instabilities (see, for example, [41]).

In the case of Fox-Goodwin, the algorithm fails regardless of the stepsize, as can be seen in Fig. 7.

This brings about an obvious conclusion: the stability issues are not caused by the non linearity of the problem itself, but by the introduction of constraints. Let us see now what happens when increasing the amplitude of the torque in the stiff pendulum problem. When doing so, the maximum angle of rotation grows and so does the rotation of the null

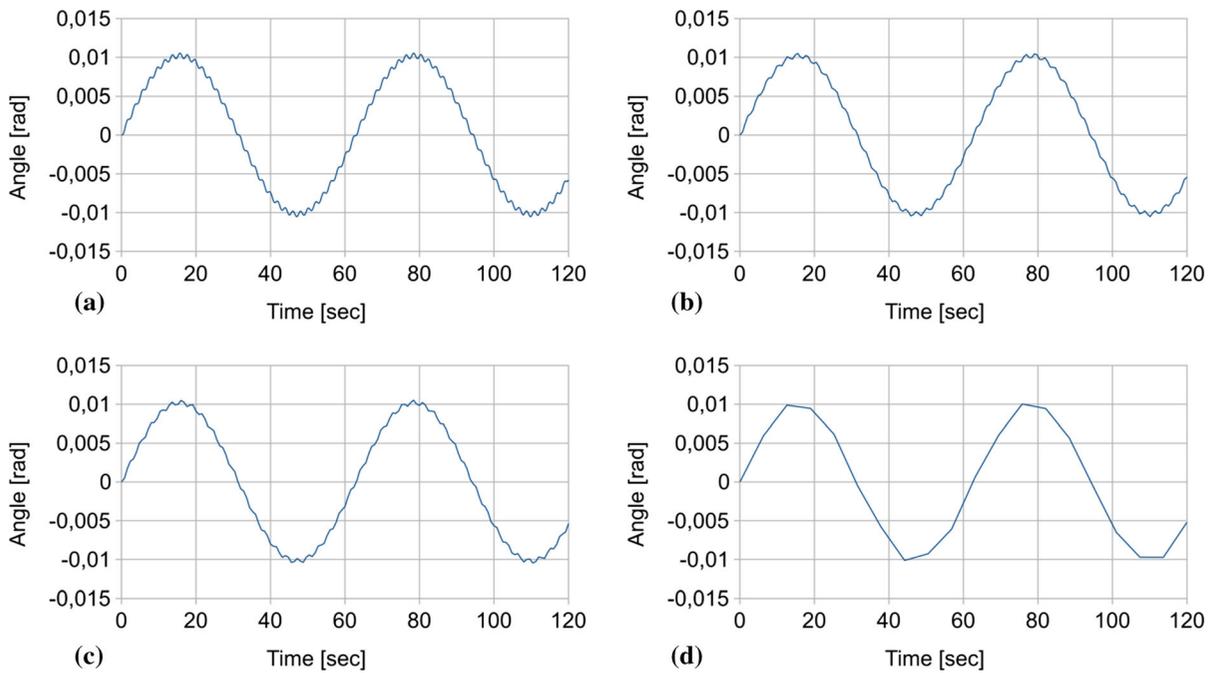

**Fig. 5** Trapezoidal rule, minimal coordinates solution for step sizes **a** 0.1 s, **b** 0.7 s, **c** 0.79 s and **d** 6 s

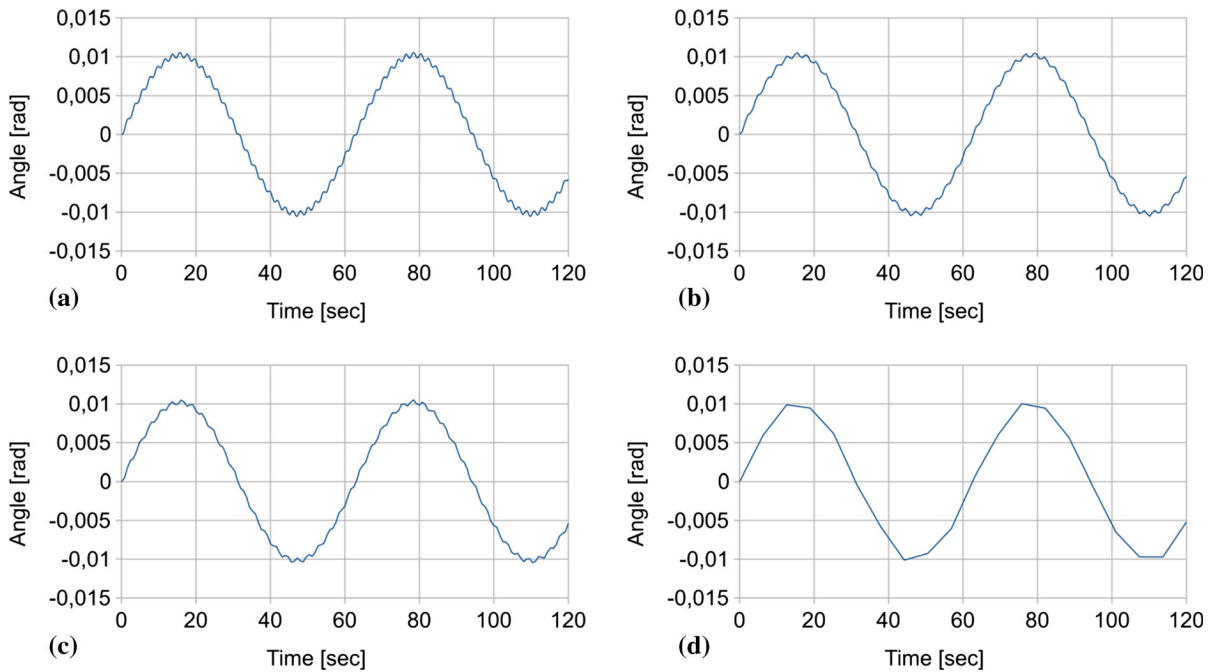

**Fig. 6** Trapezoidal rule, using constraints for step sizes **a** 0.1 s, **b** 0.7 s, **c** 0.79 s and **d** 6 s

space of the Jacobian of the constraints (see Hanke [28]). The stepsize is 0.1 s. The obtained angle for amplitudes $T_0$ of 5, 6, 8 and 9 (in Newton meter) are represented in Fig. 8:

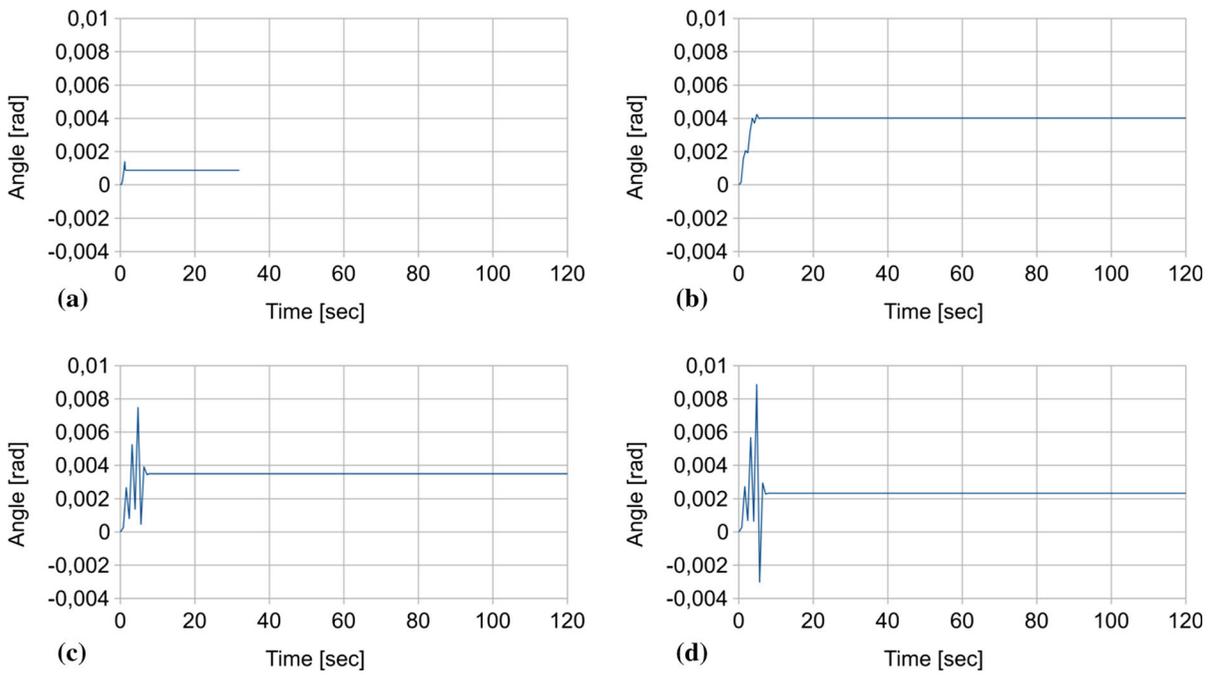

**Fig. 7** Fox-Goodwin, using constraints for step sizes **a** 0.1 s, **b** 0.6 s, **c** 0.78 s and **d** 0.79 s

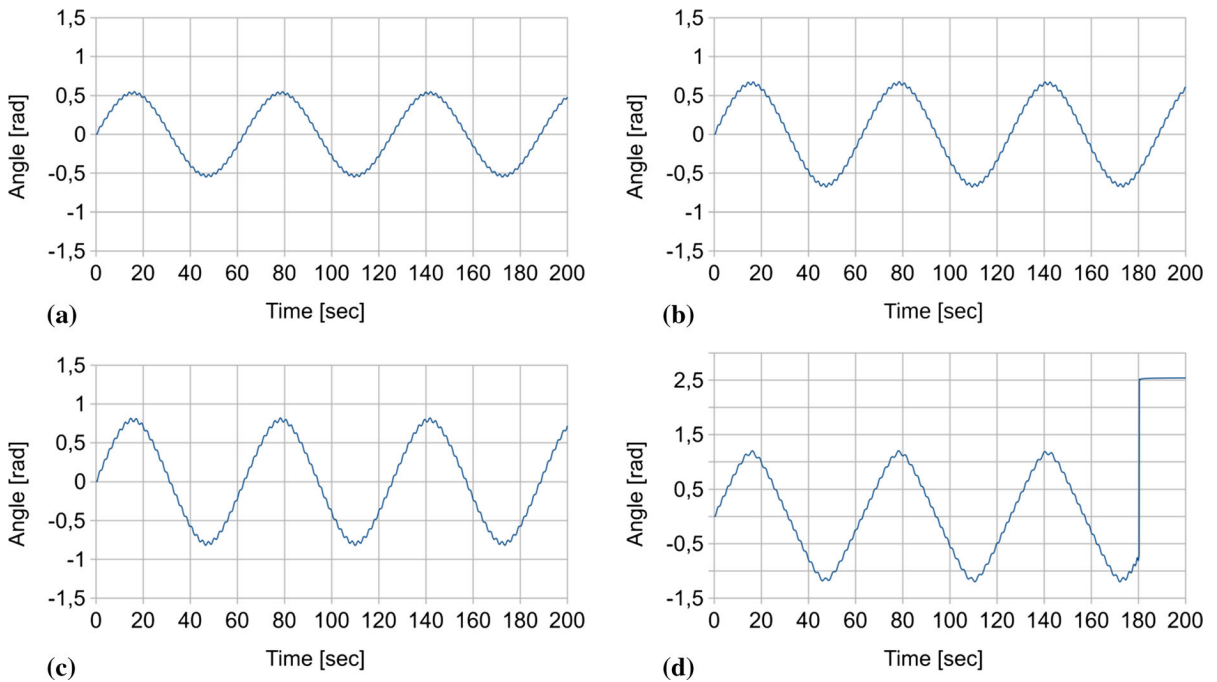

**Fig. 8** Newmark (trapezoidal rule) constrained system results for amplitudes of **a** 5Nm, **b** 6Nm, **c** 7Nm and **d** 9Nm

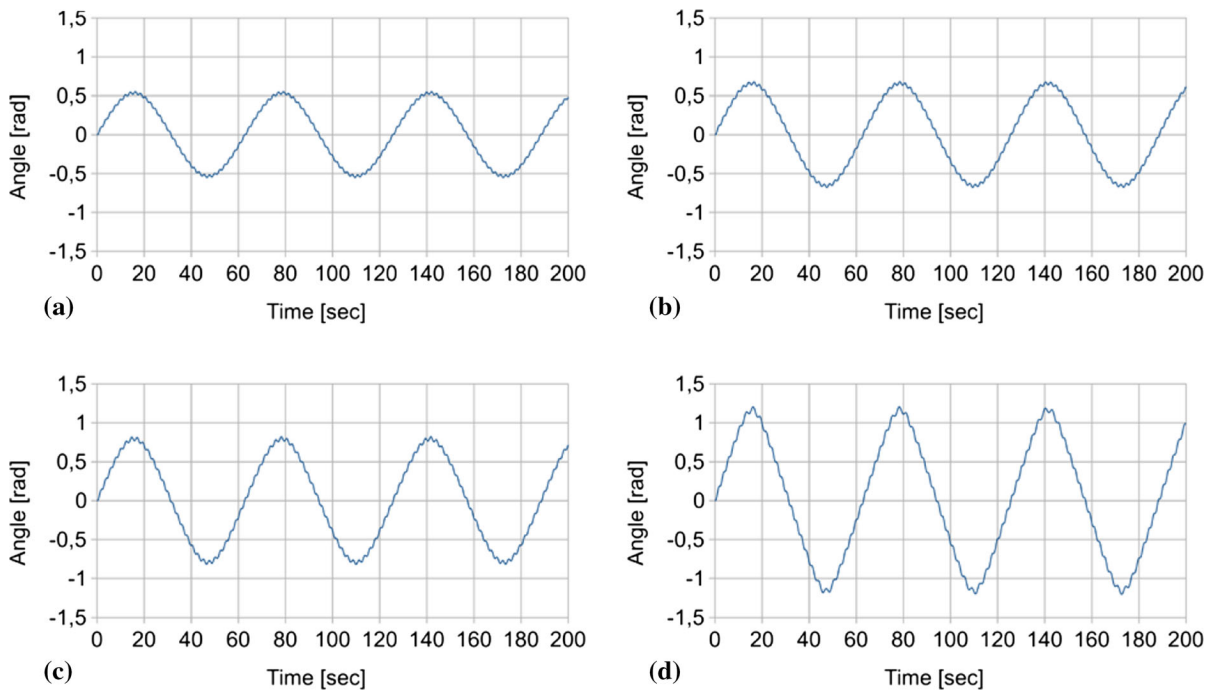

**Fig. 9** Newmark (Trapezoidal Rule) minimal coordinate results for amplitudes of **a** 5Nm, **b** 6Nm, **c** 7Nm and **d** 9Nm

The same analysis has been carried out employing the minimal coordinate approach, the same time increment (0.1 s) and Newmark parameters (Trapezoidal Rule), which yields the following results (Fig. 9):

The algorithm using constraints fails at amplitude of 9Nm. However, it does not seem to be related to instability (it rather seems to be a convergence failure). It is in fact in the reaction forces due to constraints where one can verify the instability. These forces are represented (the X component, Y is similar) in Fig. 10. In [26] the authors relate these kinds of failures of the scheme to a bad conditioning of the leading term. This issue is also approached by Hanke [28], where this leading term is defined as the null space of the Jacobian. It is also stated that its rotation compromises the stability of the scheme. Although both papers are related to index-2 DAEs, these results allow one to consider that, as the amplitude of the torque increases, so does the rotation of the pendulum and, in consequence and in the words of Hanke [28], the rotation of the nullspace of the Jacobian, being this the reason behind the failure. It is not the non linearity of the problem the reason behind the stability issues, but the introduction of constraints. In the words of Higueras, the use of minimal coordinates leads to a constant null space of the leading term, which renders the system qualified (see [26, 43]), allowing the use of a linear stability analysis to be applied. It is the introduction of constraints what makes the null space of the leading term (or, shall we say, Jacobian) to vary (or rotate) over time causing the showcased instability. When the torque amplitude is low, the rotation is quite small and thus, the problem is nearly linear. This is the reason behind the apparently correct behavior of Newmark. This looks awkward, because one could guess that, integrating in a broader space of the manifold, the average natural frequency of the system is lower and this should help stability, but the effect of the constraints prevails.

## 4 Stability considerations

For linear systems, there are several ways to study stability. Probably the most straightforward approach studies the amplification introduced by the matrix that delivers the variables in $t + \Delta t$ from those evaluated in the previous step ($t$). In this way, in the case of Newmark, one reaches, for the single degree of freedom problem:

$$\begin{Bmatrix} \Delta t^2 \ddot{x}(t+\Delta t) \\ \Delta t \dot{x}(t+\Delta t) \\ x(t+\Delta t) \end{Bmatrix} = \frac{1}{(2+4\alpha\xi\omega\Delta t + 2\beta\omega^2\Delta t^2)}$$
$$\begin{bmatrix} (4(\alpha-1)\xi\omega\Delta t + \omega^2(2\beta-1)\Delta t^2) \\ (2(1-\alpha) + 2(\beta-\alpha)\omega^2\Delta t^2) \\ ((1-2\beta) + (\alpha-2\beta)2\xi\omega\Delta t) \end{bmatrix}$$
$$\begin{matrix} (-2\omega^2\Delta t^2 - 4\xi\omega\Delta t) & (-2\omega^2\Delta t^2) \\ (2(\beta-\alpha)\omega^2\Delta t^2 + 2) & -2\alpha\omega^2\Delta t^2 \\ (2+(\alpha-\beta)4\xi\omega\Delta t) & (2+4\alpha\xi\omega) \end{matrix} \begin{Bmatrix} \Delta t^2 \ddot{x}(t) \\ \Delta t \dot{x}(t) \\ x(t) \end{Bmatrix}$$
$$= [A] \begin{Bmatrix} \Delta t^2 \ddot{x}(t) \\ \Delta t \dot{x}(t) \\ x(t) \end{Bmatrix}$$
(25)

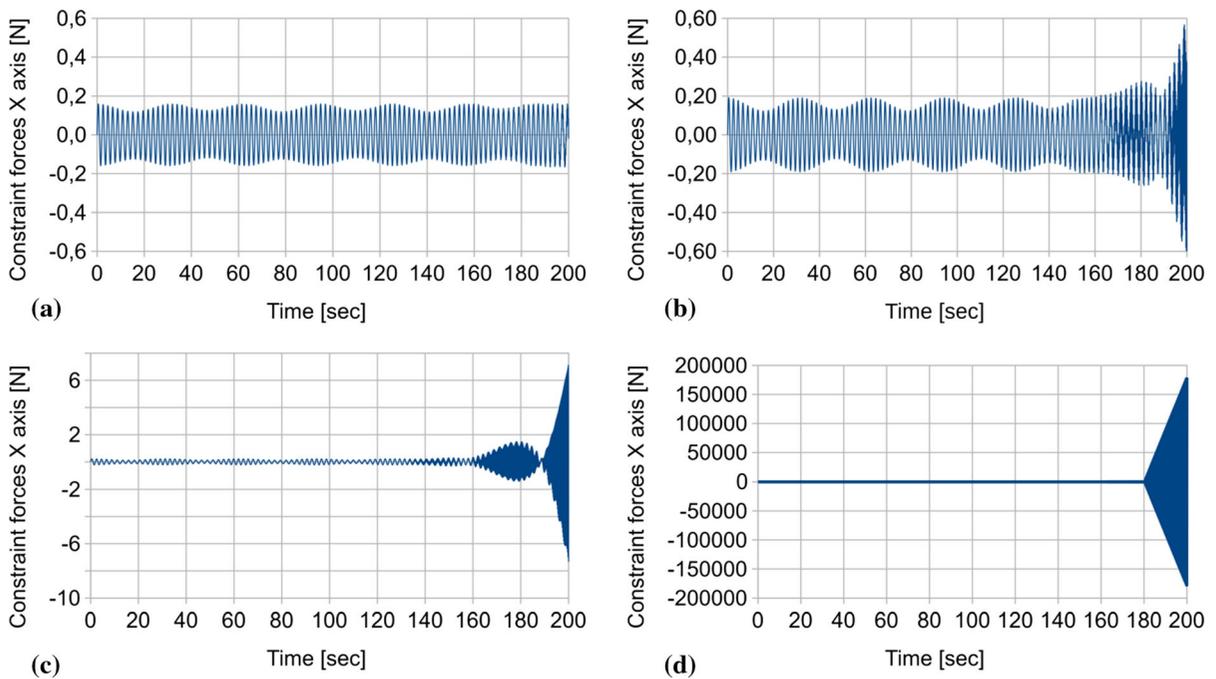

**Fig. 10** Newmark (trapezoidal rule) reaction force (X component) for amplitudes of **a** 5Nm, **b** 6Nm, **c** 7Nm and **d** 9Nm

The absolute values of the eigenvalues of matrix $[A]$ cannot be larger than 1. With this in mind, one can reach the condition expressed in the following equation:

$$\Delta t \leq \frac{1}{\omega}\sqrt{\frac{1}{\frac{1}{2}\alpha - \beta}} \quad (26)$$

This can easily be generalized to multiple degrees of freedom by decoupling the equations. This leads to Eq. (11). In a linear ODE, this equation represents the phenomena in a perfect way, because the natural frequencies of the system to integrate are constant. In nonlinear conditions, however, if one uses minimal coordinates, a violation of the equation in a single step would lead to an unstable increment of the error that might not be of importance if it does not happen repeatedly over the integration process. In any case, if this violation does not happen in any step, the system will be stable. This is also confirmed in the stiff pendulum example. This allows one to develop unconditionally stable integrators in nonlinear conditions provided that no constraints appear in the problem.

Let us now consider a constrained problem. As stated before, one can formulate constraints based on the functions to be integrated (Index 3 formulation), their first derivative (Index 2) or their second derivative (Index 1). Let us choose the most straightforward method, which is using an index 3 formulation. When trying to analyze the stability of the problem, the linearized constraints in Eq. (10) should be studied, where:

$$H_L = \frac{\partial q(x)}{\partial x} \quad (27)$$

This leads to a system that can be expressed in the following form:

$$\begin{bmatrix} M_L \\ 0 \end{bmatrix} \ddot{x}(t+\Delta t) + \begin{bmatrix} C_L \\ 0 \end{bmatrix} \dot{x}(t+\Delta t) + \begin{bmatrix} K_L \\ H_L \end{bmatrix} x(t+\Delta t) = \begin{Bmatrix} f_L + G_L^T \lambda \\ b_L \end{Bmatrix} \quad (28)$$

Now $\lambda$ can be considered as a function to be solved. Equation (28) could be written as:

$$\begin{bmatrix} M_L & 0 \\ 0 & 0 \end{bmatrix} \begin{Bmatrix} \ddot{x}(t+\Delta t) \\ \ddot{\lambda}(t+\Delta t) \end{Bmatrix} + \begin{bmatrix} C_L & 0 \\ 0 & 0 \end{bmatrix} \begin{Bmatrix} \dot{x}(t+\Delta t) \\ \dot{\lambda}(t+\Delta t) \end{Bmatrix} + \begin{bmatrix} K_L & -G_L^T \\ H_L & 0 \end{bmatrix} \begin{Bmatrix} x(t+\Delta t) \\ \lambda(t+\Delta t) \end{Bmatrix} = \begin{Bmatrix} f_L \\ b_L \end{Bmatrix} \quad (29)$$

It does not matter how small the values in $H_L$ are. The system has great probabilities of instability. If one considers Eq. (11), $\omega_{max} = \infty$, and, thus, one reaches:

$$\Delta t \leq \frac{1}{\infty} \sqrt{\frac{1}{\frac{1}{2}\alpha - \beta}} \quad (30)$$

Thus, the only way for the system to verify the stability condition is that $\frac{1}{2}\alpha - \beta \leq 0$. This is the reason behind the failure of the Fox-Goodwin approach in the pendulum example regardless of the stepsize. In the case of the classical Newmark approach ($\alpha > 0.5$ and $\beta > 0.5\alpha$), the stepsize limit is not determined and that allows the integrator to successfully solve the problem under certain conditions.

It is known that the use of formulations which introduce the constraints in a second order derivative frame tends to be more stable. This comes straightforward from the fact that, in this case, one would reach an equation in the form:

$$\begin{bmatrix} M_L & 0 \\ H_L & 0 \end{bmatrix} \begin{Bmatrix} \ddot{x}(t+\Delta t) \\ \ddot{\lambda}(t+\Delta t) \end{Bmatrix} + \begin{bmatrix} C_L & 0 \\ 0 & 0 \end{bmatrix} \begin{Bmatrix} \dot{x}(t+\Delta t) \\ \dot{\lambda}(t+\Delta t) \end{Bmatrix} + \begin{bmatrix} K_L & -G^T \\ 0 & 0 \end{bmatrix} \begin{Bmatrix} x(t+\Delta t) \\ \lambda(t+\Delta t) \end{Bmatrix} = \begin{Bmatrix} f_L \\ c_L \end{Bmatrix} \quad (31)$$

This introduces nonzeroes in the mass matrix instead of introducing them in the stiffness matrix. The problem here lies in the drift-off of the function to be integrated, which must be somehow corrected. One can introduce Baumgarte stabilization, or projection, or any method to avoid the drift but, anyway, these methods introduce a relation among the functions in $t + \Delta t$ and those in $t$, which, in any case will numerically behave as a stiffness, rendering again the stability analysis useless. This is a major concern, because in these cases one cannot predict stability in a proper way. In addition, another problem of these stabilization methods is that they introduce more equations than variables in each integration step, which drives to the need of finding a solution of compromise.

Note that the mentioned impossibility of performing a proper stability analysis and the associated issues appear, as Eq. (31) shows, in implicit methods; explicit methods, on the other hand, do not present this problem but, unfortunately, there is no unconditionally stable explicit method. For example, in the central differences method the constraints do not appear as a stiffness and, thus, an integrator based on central differences has predictable conditional stability.

One could also argue that Eqs. (29) and (31) could face numerical issues. In the usual approaches that are used to solve them, these are avoided, as mentioned by Nocedal [11].

## 5 A new approach for solving the problem

The unexpected stability issues seem to arise not due to the nonlinearity of the problem, but because of the introduction of nonlinear constraints. Problems solved using minimal coordinates do not exhibit unpredictable behavior, while those solved with non-linear constraints do. A quite interesting study in the subject can be found in [26]. In that paper, the use of minimal coordinates would be considered to lead to a qualified DAE [26, 43]. One could also state that, in minimal coordinates, no constraints are to be dealt with, and thus, the problem is an ODE. However, one cannot always rely on minimal coordinates, because this usually involves formulations that are not general. The question arisen is: is there a way to integrate a DAE avoiding this problem? It is obvious that not in the usual form. Instead of integrating the variables in which the DAE is formulated, an alternative is to integrate in the manifold (which is, in fact, what minimal coordinates do). To accomplish this procedure in a DAE, one can integrate in a linearization of the manifold performed in the step to integrate. These ideas can be applied to most implicit integrators, but, for the sake of clarity, we will use a Newmark approach, using analytical derivatives. A similar approach was presented in [34], but applied to an explicit method. Due to the fact that explicit methods do not exhibit the stability behavior that implicit methods do, the approach in [34] cannot be used here.

Let us consider the DAE formulated with eqns. (3) and (2). Here the full restriction treatment will be presented. A Taylor series linearization of (2) leads to:

$$\boldsymbol{q}(\boldsymbol{x}) \approx \boldsymbol{q}(\boldsymbol{x}_0) + \left.\frac{\partial \boldsymbol{q}(\boldsymbol{x})}{\partial \boldsymbol{x}}\right|_{\boldsymbol{x}_0} (\boldsymbol{x} - \boldsymbol{x}_0) = 0 \quad (32)$$

Taking into account (27), one can write:

$$\boldsymbol{H}_L(\boldsymbol{x}_0)\boldsymbol{x} = \boldsymbol{H}_L(\boldsymbol{x}_0)\boldsymbol{x}_0 - \boldsymbol{q}(\boldsymbol{x}_0) = \boldsymbol{b}_L \quad (33)$$

If one considers this equation for obtaining $\boldsymbol{x}(t + \Delta t)$, $\boldsymbol{x}_0$ will be the current estimation of $\boldsymbol{x}(t + \Delta t)$. With this in mind, one can write:

$$\boldsymbol{x}(t + \Delta t) = \boldsymbol{x}_p + \boldsymbol{N}_H \boldsymbol{\alpha}(t + \Delta t) \quad (34)$$

where $\boldsymbol{N}_H$ is the null space of $\boldsymbol{H}_L$, and $\boldsymbol{x}_p$ is a particular solution of Eq. (10). Any expression of $\boldsymbol{N}_H$ and $\boldsymbol{x}_p$ would theoretically lead to the correct solution. An obvious choice for $\boldsymbol{x}_p$ is the minimal least squares solution. This choice will always provide lead to good numerical conditioning. One could argue that the use of a sparse solution for $\boldsymbol{x}_p$ could be of use, but being $\boldsymbol{x}_p$ a vector, little advantage can be obtained from this. The use of different expressions for $\boldsymbol{x}_p$ should generate proper results, provided that the elements of $\boldsymbol{x}_p$ are kept small. For $\boldsymbol{N}_H$ the use of a fundamental base will bring about a sparse expression, which will reduce computational cost, as mentioned by Coleman [44]. The important fact is for $\boldsymbol{\alpha}(t + \Delta t)$ to be a set of minimal coordinates which are only valid for the step currently being integrated. The hereby presented method is based on integrating on these coordinates, which change from step to step. In the case of both scleronomic and holonomic constraints, one can write:

$$\dot{\boldsymbol{q}}(\boldsymbol{x}, \dot{\boldsymbol{x}}) \simeq \dot{\boldsymbol{q}}(\boldsymbol{x}_0, \dot{\boldsymbol{x}}_0) + \left.\frac{\partial \dot{\boldsymbol{q}}}{\partial \boldsymbol{x}}\right|_{\boldsymbol{x}_0, \dot{\boldsymbol{x}}_0} (\boldsymbol{x} - \boldsymbol{x}_0)$$
$$+ \left.\frac{\partial \boldsymbol{q}}{\partial \boldsymbol{x}}\right|_{\boldsymbol{x}_0} (\dot{\boldsymbol{x}} - \dot{\boldsymbol{x}}_0) = 0 \quad (35)$$

This can be written in the form:

$$\boldsymbol{H}(\boldsymbol{x}_0)\dot{\boldsymbol{x}} = -\dot{\boldsymbol{q}}(\boldsymbol{x}_0, \dot{\boldsymbol{x}}_0) + \boldsymbol{H}(\boldsymbol{x}_0)\dot{\boldsymbol{x}}_0 + \dot{\boldsymbol{H}}(\boldsymbol{x}_0, \dot{\boldsymbol{x}}_0)\boldsymbol{x}_0$$
$$- \dot{\boldsymbol{H}}(\boldsymbol{x}_0, \dot{\boldsymbol{x}}_0)\boldsymbol{x} \quad (36)$$

Now one introduces (34), reaching:

$$\boldsymbol{H}(\boldsymbol{x}_0)\dot{\boldsymbol{x}} = -\dot{\boldsymbol{q}}(\boldsymbol{x}_0, \dot{\boldsymbol{x}}_0) + \boldsymbol{H}(\boldsymbol{x}_0)\dot{\boldsymbol{x}}_0 + \dot{\boldsymbol{H}}(\boldsymbol{x}_0, \dot{\boldsymbol{x}}_0)\boldsymbol{x}_0$$
$$- \dot{\boldsymbol{H}}(\boldsymbol{x}_0, \dot{\boldsymbol{x}}_0)\boldsymbol{x}_p - \dot{\boldsymbol{H}}(\boldsymbol{x}_0, \dot{\boldsymbol{x}}_0)\boldsymbol{N}_H \boldsymbol{\alpha} \quad (37)$$

The system (37) can be solved leading to:

$$\dot{\boldsymbol{x}} = \dot{\boldsymbol{x}}_P + \boldsymbol{N}_H \dot{\boldsymbol{\alpha}} + \dot{\boldsymbol{X}}_P \boldsymbol{\alpha} \quad (38)$$

Being $\dot{\boldsymbol{x}}_p$ a solution to:

$$\begin{aligned}\boldsymbol{H}(\boldsymbol{x}_0)\dot{\boldsymbol{x}} &= -\dot{\boldsymbol{q}}_p(\boldsymbol{x}_0, \dot{\boldsymbol{x}}_0) + \boldsymbol{H}(\boldsymbol{x}_0)\dot{\boldsymbol{x}}_0 + \dot{\boldsymbol{H}}(\boldsymbol{x}_0, \dot{\boldsymbol{x}}_0)\boldsymbol{x}_0 \\ &\quad - \dot{\boldsymbol{H}}(\boldsymbol{x}_0, \dot{\boldsymbol{x}}_0)\boldsymbol{x}_p \\ &= -\dot{\boldsymbol{q}}_p(\boldsymbol{x}_0, \dot{\boldsymbol{x}}_0) + \boldsymbol{H}(\boldsymbol{x}_0)\dot{\boldsymbol{x}}_0 \\ &\quad + \dot{\boldsymbol{H}}(\boldsymbol{x}_0, \dot{\boldsymbol{x}}_0)(\boldsymbol{x}_0 - \boldsymbol{x}_p)\end{aligned} \quad (39)$$

And $\dot{\boldsymbol{X}}_P$ a solution to:

$$\boldsymbol{H}(\boldsymbol{x}_0)\dot{\boldsymbol{X}} = -\dot{\boldsymbol{H}}(\boldsymbol{x}_0, \dot{\boldsymbol{x}}_0)\boldsymbol{N}_H \quad (40)$$

Here it must be pointed out that neither $\dot{\boldsymbol{x}}_p$ nor $\dot{\boldsymbol{X}}_P$ are obtained as derivatives of any function, but as solutions of linear systems. The use of the dot has been decided in order to keep unit coherence. Now one needs to tackle accelerations. By taking again a derivative:

$$\ddot{\boldsymbol{q}}(\boldsymbol{x}, \dot{\boldsymbol{x}}, \ddot{\boldsymbol{x}}) \simeq \ddot{\boldsymbol{q}}(\boldsymbol{x}_0, \dot{\boldsymbol{x}}_0, \ddot{\boldsymbol{x}}_0) + \left.\frac{\partial \ddot{\boldsymbol{q}}}{\partial \boldsymbol{x}}\right|_{\boldsymbol{x}_0, \dot{\boldsymbol{x}}_0, \ddot{\boldsymbol{x}}_0} (\boldsymbol{x} - \boldsymbol{x}_0)$$
$$+ 2\left.\frac{\partial \dot{\boldsymbol{q}}}{\partial \boldsymbol{x}}\right|_{\boldsymbol{x}_0, \dot{\boldsymbol{x}}_0} (\dot{\boldsymbol{x}} - \dot{\boldsymbol{x}}_0) + \left.\frac{\partial \boldsymbol{q}}{\partial \boldsymbol{x}}\right|_{\boldsymbol{x}_0} (\ddot{\boldsymbol{x}} - \ddot{\boldsymbol{x}}_0)$$
$$= 0 \quad (41)$$

This can be rewritten as:

$$\ddot{q}(x_0, \dot{x}_0, \ddot{x}_0) + \ddot{H}(x_0, \dot{x}_0, \ddot{x}_0)(x - x_0)$$
$$+ 2\dot{H}(x_0, \dot{x}_0)(\dot{x} - \dot{x}_0) + H(x_0)(\ddot{x} - \ddot{x}_0)$$
$$= 0 \qquad (42)$$

Reordering:

$$H(x_0)\ddot{x} = -\ddot{q}(x_0, \dot{x}_0, \ddot{x}_0) + \ddot{H}(x_0, \dot{x}_0, \ddot{x}_0)x_0$$
$$+ 2\dot{H}(x_0, \dot{x}_0)\dot{x}_0 + H(x_0)\ddot{x}_0$$
$$- \ddot{H}(x_0, \dot{x}_0, \ddot{x}_0)x - 2\dot{H}(x_0, \dot{x}_0)\dot{x} \qquad (43)$$

Now one can substitute eqns. (34) and (38):

$$H(x_0)\ddot{x} = -\ddot{q}(x_0, \dot{x}_0, \ddot{x}_0) + \ddot{H}x_0 + 2\dot{H}\dot{x}_0 + H\ddot{x}_0$$
$$- \ddot{H}x_p - 2\dot{H}\dot{x}_p - 2\dot{H}N_H\dot{\alpha}$$
$$+ \left(-2\dot{H}\dot{X}_p - \ddot{H}N_H\right)\alpha \qquad (44)$$

So one can solve:

$$\ddot{x} = \ddot{x}_p + \ddot{X}_{p1}\alpha + X_{p2}\dot{\alpha} + N_H\ddot{\alpha} \qquad (45)$$

Being $\ddot{x}_p$ a particular solution to:

$$H(x_0)\ddot{x} = -\ddot{q}(x_0, \dot{x}_0, \ddot{x}_0) + \ddot{H}(x_0, \dot{x}_0, \ddot{x}_0)(x_0 - x_p)$$
$$+ 2\dot{H}(x_0, \dot{x}_0)(\dot{x}_0 - \dot{x}_p) + H(x_0)\ddot{x}_0 \qquad (46)$$

In turn, $\ddot{X}_{p1}$ is a solution to:

$$H(x_0)\ddot{X} = -2\dot{H}(x_0, \dot{x}_0)\dot{X}_p - \ddot{H}(x_0, \dot{x}_0, \ddot{x}_0)N_H \qquad (47)$$

And $\ddot{X}_{p2}$ is a solution to:

$$H(x_0)\ddot{X} = -2\dot{H}(x_0, \dot{x}_0)N_H \qquad (48)$$

One can notice:

$$\dot{X}_{P2} = 2\dot{X}_P \qquad (49)$$

Regarding $x_p$, $\dot{X}_P$, $\dot{x}_p$, $\ddot{X}_{p1}$ and $\ddot{x}_p$, it is important to state that, although any solution will suffice, it is numerically important not to let these solutions take any value (it could be too large and cause numerical issues), thus, it is recommended to use minimal norm values. In the event of redundant constraints, it is recommended to use minimal least squares values, even considering that they should be compatible.

Now one must formulate the equilibrium equation in terms of $\alpha$. From Eq. (6), and introducing eqs. (34), (38) and (42), one can write:

$$M_L(\ddot{x}_p + \ddot{X}_{p1}\alpha + X_{p2}\dot{\alpha} + N_H\ddot{\alpha})$$
$$+ C_L(\dot{x}_p + N_H\dot{\alpha} + \dot{X}_p\alpha) + K_L(x_p + N_H\alpha)$$
$$= f_L + G_L^T\lambda \qquad (50)$$

Reordering:

$$M_L N_H \ddot{\alpha} + \left(C_L N_H + 2M_L \dot{X}_p\right)\dot{\alpha}$$
$$+ \left(K_L N_H + C_L \dot{X}_p + M_L \ddot{X}_{p1}\right)\alpha$$
$$= f_L - K_L x_p - C_L \dot{x}_p - M_L \ddot{x}_p + G_L^T\lambda \qquad (51)$$

This is the system to be actually integrated. This system is not constrained, because any value given to $\alpha(t + \Delta t)$, $\dot{\alpha}(t + \Delta t)$ and $\ddot{\alpha}(t + \Delta t)$ will lead to values of $x(t + \Delta t)$, $\dot{x}(t + \Delta t)$ and $\ddot{x}(t + \Delta t)$ satisfying the linearized constraints. Thus, in an iterative approach they will converge to satisfying the original (not linearized) constraints. It is important to remark that the approach will not only converge to the verification of either the constraints expressed in terms of the function, its first derivative or its second derivative, but also all of them together. Once computed $\alpha(t + \Delta t)$, $\dot{\alpha}(t + \Delta t)$ and $\ddot{\alpha}(t + \Delta t)$, one can use eqs. (34), (38) and (42) to obtain $x(t + \Delta t)$, $\dot{x}(t + \Delta t)$ and $\ddot{x}(t + \Delta t)$. The process will be restarted with these initial values until convergence is achieved, leading to the next timestep. To avoid ill-conditioned systems, one can use the null space of $G$ to reduce the system:

$$N_G^T M_L N_H \ddot{\alpha} + \left(N_G^T C_L N_H + 2N_G^T M_L \dot{X}_p\right)\dot{\alpha}$$
$$+ \left(N_G^T K_L N_H + N_G^T C_L \dot{X}_p + N_G^T M_L \ddot{X}_{p1}\right)\alpha$$
$$= N_G^T f_L - N_G^T K_L x_p - N_G^T C_L \dot{x}_p - N_G^T M_L \ddot{x}_p$$
$$+ N_G^T G_L^T\lambda \qquad (52)$$

A commentary is due here. The novelty of the hereby presented method is not the use of the null space of $G_L$ to reduce the system, but on the fact that $\alpha$ and its derivatives are what is being integrated instead of $x$ and its derivatives, as can be seen in Eq. (52). This differs from the usual way, in which first the integration is applied to an equation in the form of (23) effectively integrating $x$, and then constraints of either position, velocity or acceleration are applied. Therefore, in the usual way only one kind of constraints is applied and fulfilled.

Another important issue is that the frequencies to be taken into account for stability analysis are those of the system to be actually integrated, and, thus, those should be obtained from the reduced mass matrix $M_R = N_G^T M_L N_H$ and the reduced stiffness matrix $K_R = N_G^T K_L N_H + N_G^T C_L \dot{X}_p + N_G^T M_L \ddot{X}_{p1}$.

The computational cost can be increased because one needs to obtain $N_G^T K_L N_H + N_G^T C_L \dot{X}_p + N_G^T M_L \ddot{X}_{p1}$ and $N_G^T C_L N_H + 2N_G^T M_L \dot{X}_p$. But $N_G^T K_L N_H$, $N_G^T C_L N_H$, $N_G^T M_L$ and $N_G^T C_L$ have to be obtained anyway, one only needs to compute three products of reduced matrices. The cost of obtaining $\ddot{X}_{p1}$ and $\dot{X}_p$ is also noticeable, but one must take into account that the factorization of $N_H$ is also needed regardless of the method, so it introduces no additional cost. Thus, the method does not change the order of the computational cost.

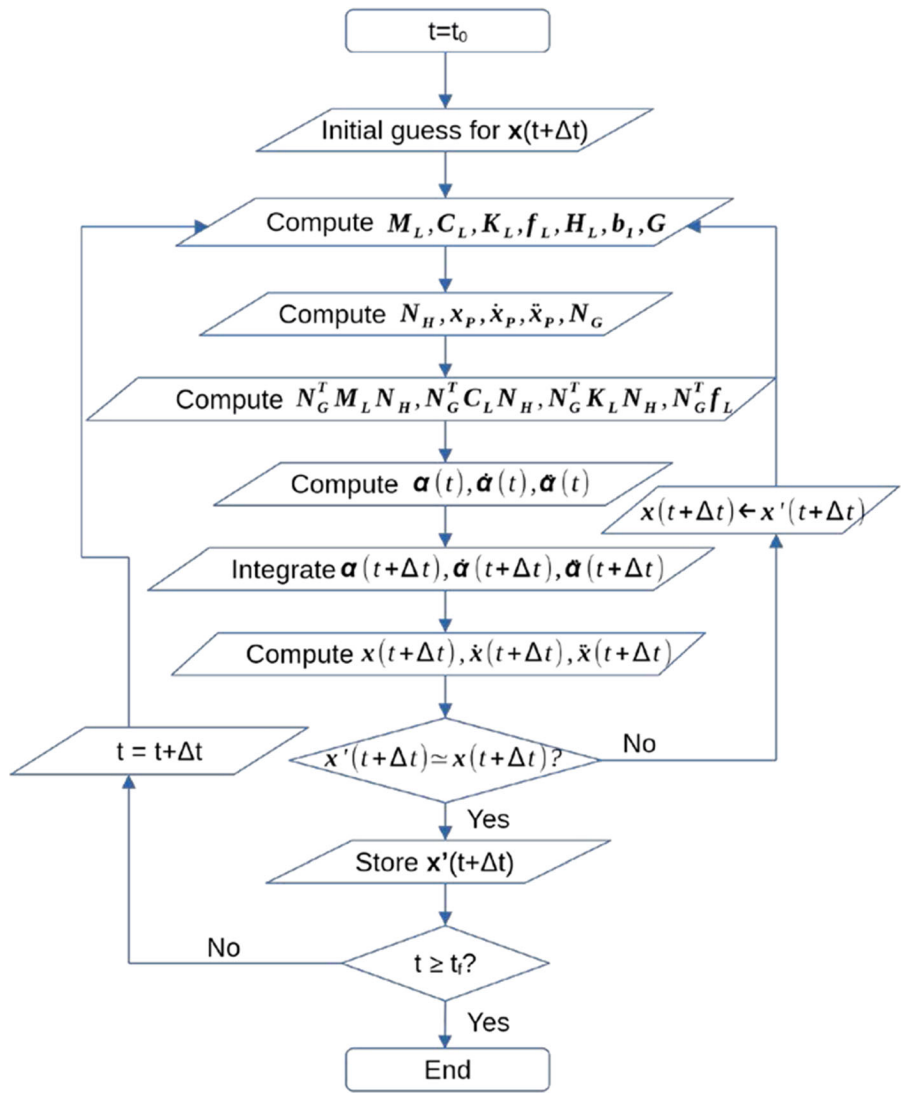

**Fig. 11** Flowchart of the proposed algorithm

Lastly, the term $N_G^T G_L^T \lambda$ can be fully eliminated in explicit methods due to the orthogonality property. Implicit methods, however, require a previous linearization of the term $G_L^T \lambda$, and only part of this term will be eliminated due to the mentioned orthogonality property. The remaining terms of the linearization will be included as part of the stiffness matrix $K_L$ and $f_L$.

One of the key points of the approach presented here is to notice that, if $\alpha$ is a variable to be integrated, $\alpha(t)$, $\dot{\alpha}(t)$ and $\ddot{\alpha}(t)$ are also needed. It could be tempting to use the values obtained in the previous step, but this cannot be done, because the linearization is not the same from step to step. Furthermore, that would not be integrating in the linearization. In order to obtain $\alpha(t)$, $\dot{\alpha}(t)$ and $\ddot{\alpha}(t)$ one generalizes the eqs. (34), (38) and (42), which were formulated in $t + \Delta t$ to the following eqns.:

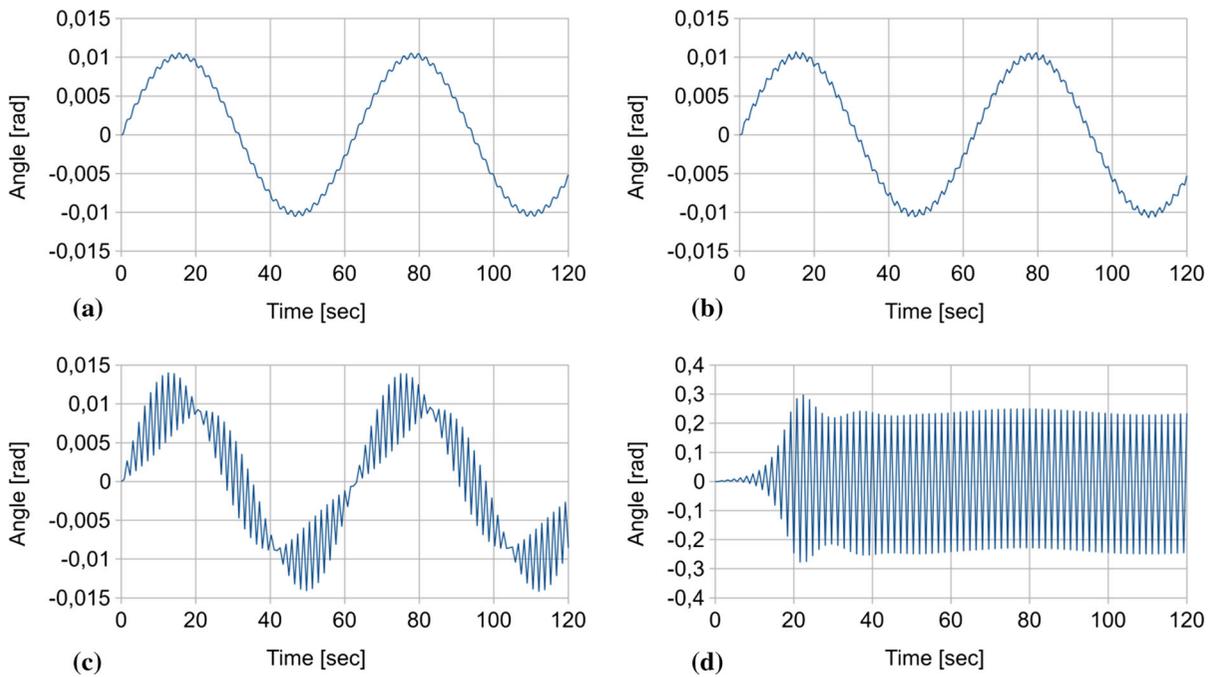

**Fig. 12** Fox-Goodwin, new formulation, for stepsizes **a** 0.1 s, **b** 0.6 s, **c** 0.78 s and **d** 0.79 s

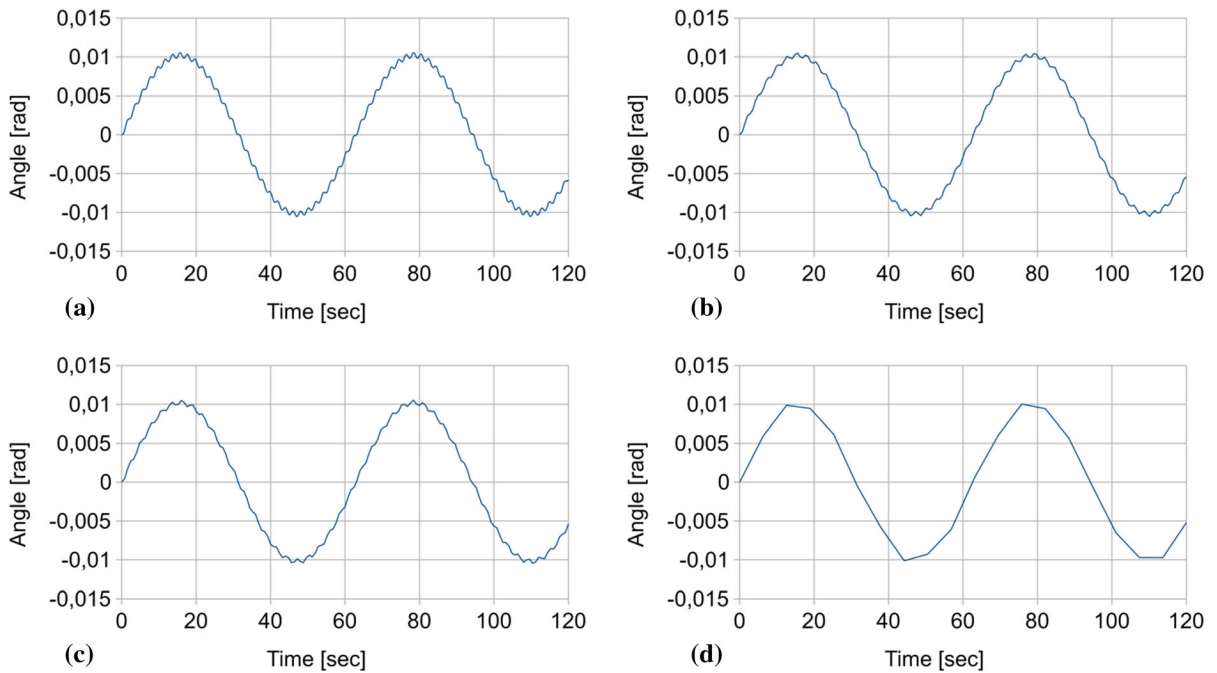

**Fig. 13** Trapezoidal rule, using constraints for stepsizes **a** 0.1 s, **b** 0.7 s, **c** 0.79 s and **d** 6 s

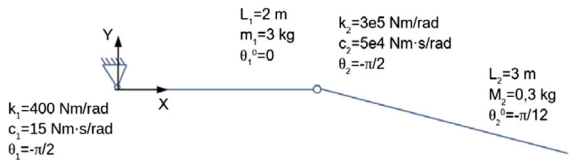

**Fig. 14** Double pendulum

$$x(t) = x_p + N_H \alpha(t) \tag{53}$$

$$\dot{x}(t) = \dot{x}_P + N_H \dot{\alpha}(t) + \dot{X}_p \alpha(t) \tag{54}$$

$$\ddot{x}(t) = \ddot{x}_p + \ddot{X}_{p1} \alpha(t) + \dot{X}_{p2} \dot{\alpha}(t) + N_H \ddot{\alpha}(t) \tag{55}$$

As the unknowns are $\alpha(t)$, $\dot{\alpha}(t)$ and $\ddot{\alpha}(t)$, it can be written as follows:

$$N_H \alpha(t) = x(t) - x_p \tag{56}$$

$$N_H \dot{\alpha}(t) = \dot{x}(t) - \dot{x}_P - \dot{X}_p \alpha(t) \tag{57}$$

$$N_H \ddot{\alpha}(t) = \ddot{x}(t) - \ddot{x}_P - \ddot{X}_{p1} \alpha(t) - \dot{X}_{p2} \dot{\alpha}(t) \tag{58}$$

The eqs. (56), (57) and (58) are overdetermined and usually incompatible, because the linearization in $t + \Delta t$ will not usually satisfy the nonlinear constraints in t. Thus, the least squares solutions for these equations can be computed, which is equivalent to project the coordinates in t to the linearization in $t + \Delta t$. Obviously, if the function to integrate is continuous, the smaller $\Delta t$ is, the better the projection fits to the original function. As a consequence, convergence of this operation seems unquestionable. With all this in mind, the flowchart of the algorithm is represented in Fig. 11.

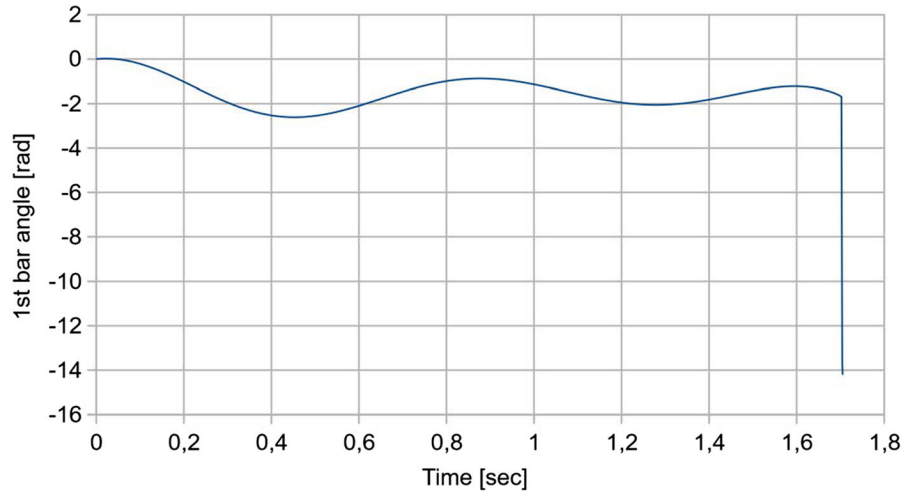

**Fig. 15** Results of Newmark integration based on displacements (angle of the 1st truss)

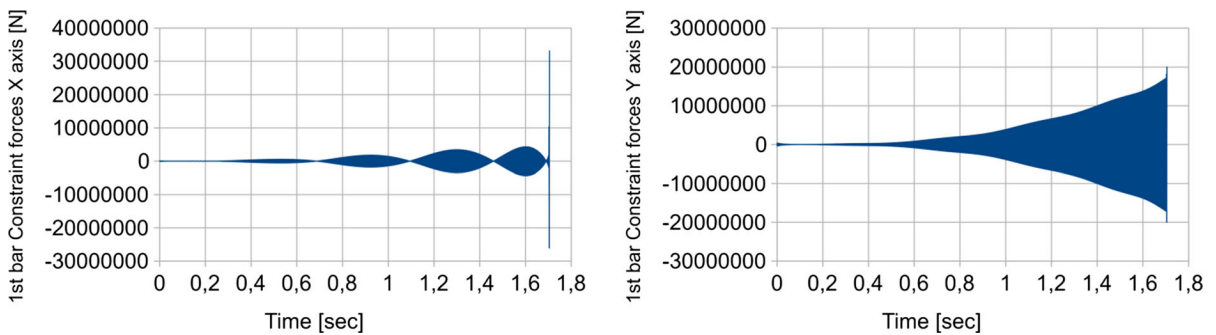

**Fig. 16** Value of the constraint forces in the fixed joint. Left: in the X axis. Right: in the Y axis

With this approach, $\alpha$ is a set of minimal coordinates to be integrated. Its physical interpretation is linked to the null space of $H$, that is, $N_H$. The product $N_H\alpha$ (and also $N_H\dot\alpha$ and $N_H\ddot\alpha$) represents the direction which is tangent to the displacement (or velocity or acceleration fields, correspondingly) of the system in that instant.

## 6 Numerical experiments

For these numerical experiments, the Octave environment has been used. This will not lead to a proper measure of the efficiency of the method, provided that Octave does not use even a JIT compiler to improve the iterative processes. The use of Matlab would considerably improve on this, but would still not lead to a proper assessment of the efficiency.

The first example is obviously the simple stiff pendulum discussed before. The results achieved, using Fox-Goodwin parameters, for stepsizes of 0.1 s, 0.6 s, 0.78 s, 0.79 s are presented in Fig. 12:

If one compares these results to those presented in Fig. 4 and 7, one can see that the new method, being general coordinate based, behaves like a minimal coordinate approach in terms of stability. It is important to remark that the classical approach presented fails even in the smallest stepsize, as shown in Fig. 7.

The trapezoidal rule also behaves in a correct way, as expected (Fig. 13).

The second example in this document is the double pendulum analyzed by Negrut [35]. In the paper the authors had to use quite conservative values of the parameters: $\alpha = 0.75$ and $\beta = 0.390625$, which delivers a first order algorithm. Even under these conditions, the authors found the method quite competitive against BDF integrators. The used stepsize is $1e-3s$. Although the authors do not explicitly state it, the reason of such a conservative set of parameters seems related to stability. The system is depicted in Fig. 14. Both joints include torsional spring-dampers. The system is subjected to gravity in the negative y axis. The considered acceleration is $g = 9.81 m/s^2$. Total integration time is 10 s.

First, the classical approach along with a displacement-based formulation is tested. The values used are $\alpha = 0.5000005$ and $\beta = 0.2500006$, which are clearly in the unconditionally stable region. The obtained results are presented in Fig. 15:

What is interesting from this example are the values of the constraint forces in the fixed joint (see Fig. 16). They show clearly the effect of instability and seem to agree with the hypothesis presented before. That is, the constraints are the reason inducing the instability, not the nonlinearity of the equilibrium equation.

The values of $\alpha$ and $\beta$ used in this example are greater than the bare minimum required for Newmark to be unconditionally stable just to avoid any doubt about issues related to discrete mathematics. The use of $\alpha = 0.5$ and $\beta = 0.25$ also brings to failure.

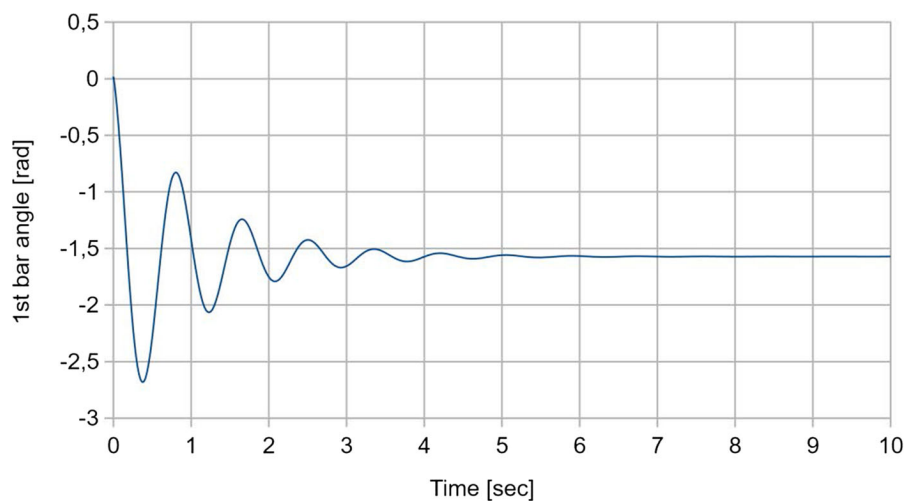

**Fig. 17** Results using the new approach: FG with $5e-4s$ stepsize

Depending on the particular implementation of the algorithm (based in acceleration with stabilization or so) the results might vary but, in any case, with the classical approach there is no way of predicting a priori the stability of the method. Using constraints formulated in terms of the positions (the function to be integrated) along with a choice of $\alpha = 0,5001$ and $\beta = 0,2507$, and a stepsize of $1e-3s$, the problem is successfully solved. But it is not quite comfortable to rely on a try-error scheme to determine a working set of parameters.

However, using the new approach the results change dramatically. One can even use a Fox-Goodwin (FG) approach. Using $\alpha = 0.5$ and $\beta = 1/12$, with a stepsize of $5e-4s$, one reaches the result shown in Fig. 17.

An additional advantage of the method is the fact that constraints are numerically fulfilled for position, velocity and acceleration (see Fig. 18, 19 and 20). The euclidean norm of the constraint errors are always below $3e-14$ for position and velocity constraints, whereas for acceleration it is under $1e-10$. This high value happens only at the beginning of the simulation, and seems related to the high accelerations produced due to initial conditions.

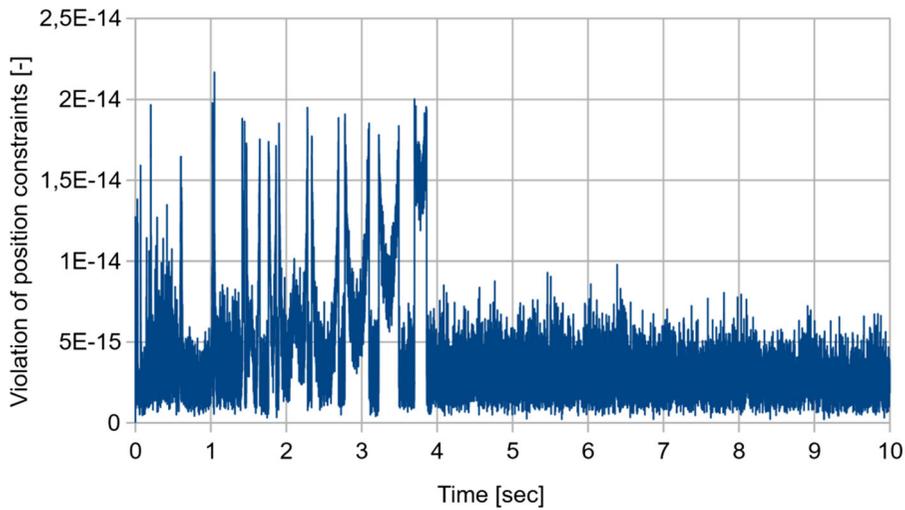

**Fig. 18** Evolution of the norm of the position constraint error

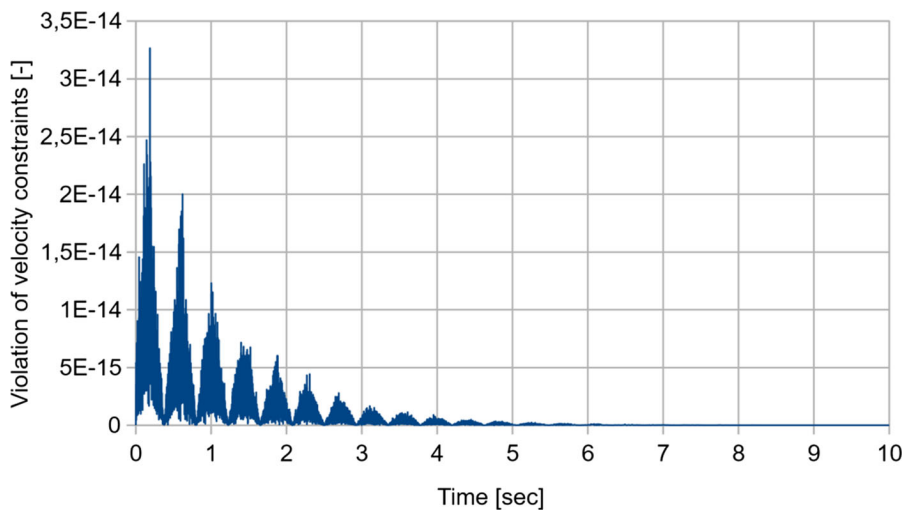

**Fig. 19** Evolution of the norm of the velocity constraint error

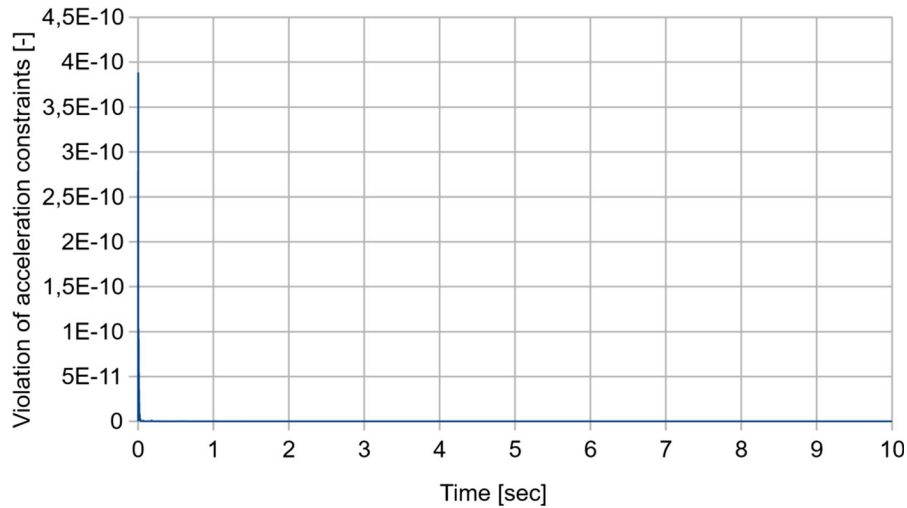

**Fig. 20** Evolution of the norm of the acceleration constraint error

If one computes the maximum natural frequencies of the system the graph in Fig. 21 is obtained.

The high value of the natural frequencies at the beginning is due to the initial values. In any case, one can compute from the linear stability analysis that a $5e-3s$ stepsize can be used along with FG. This is five times the timestep used in [35]. The results with this configuration are presented in Fig. 22.

Clearly, the result is not numerically good, mainly due to the fact that at the beginning of the simulation, there is a large variation in the accelerations which require a considerably smaller timestep. The interesting point here is that, yet, the method is stable in these conditions.

The last example to be presented is the Andrews Squeezer mechanism, which can be found in [45] and a good description with results can be spotted in the IFTOMM Multibody Benchmark.

Coordinates for the fixed joints A, B, C and O are shown in Table 1.

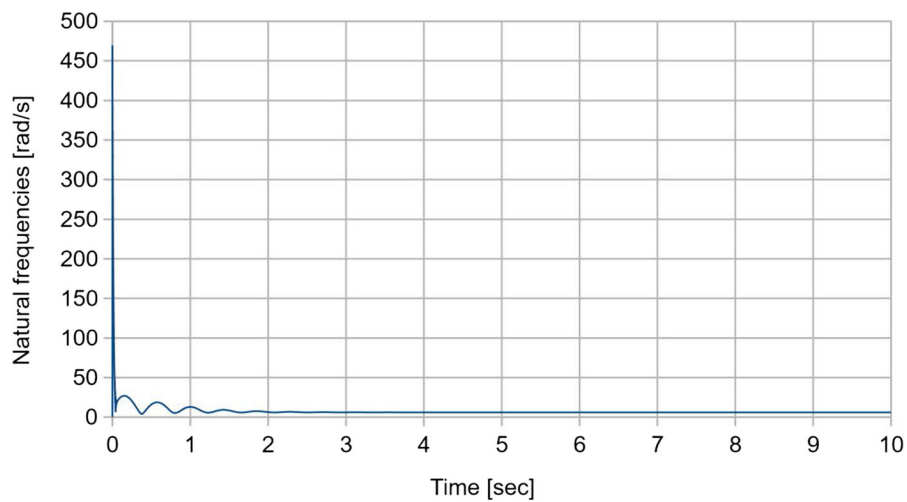

**Fig. 21** Evolution of the maximum natural frequencies of the double pendulum

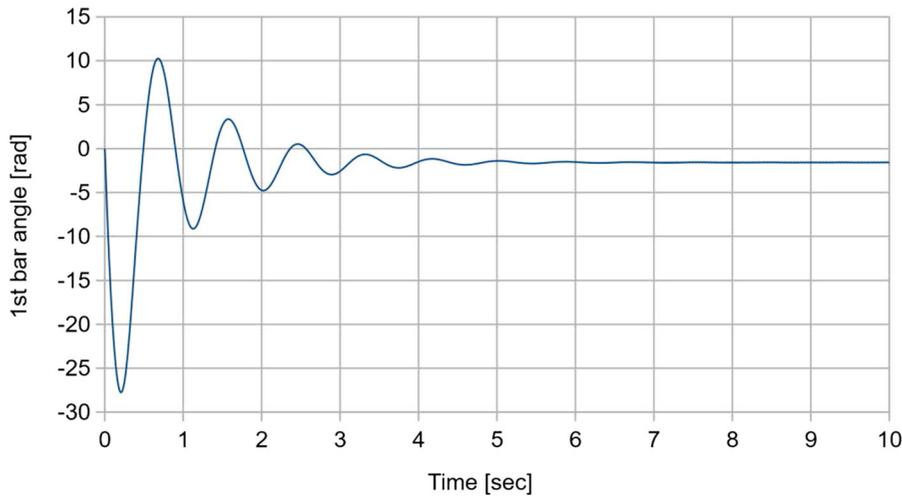

**Fig. 22** Results using the new approach: FG with $5e-3s$ stepsize

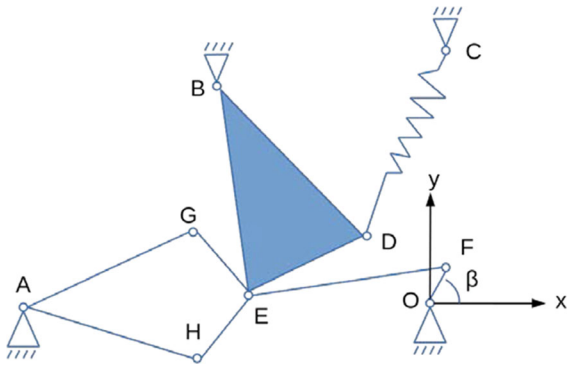

**Fig. 23** Andrews Squeezer Mechanism

Trusses have the following mechanical properties:

The coordinates for the centers of gravity of the trusses are indicated in local coordinates along the line (Table 2).

The solid BDE has a mass of 0.02373 kg and 5.255e-6 kg·m$^2$ of inertia. BE is 0.035 m length. Spring stiffness is $k = 4530 N/m$, and its undeformed length is $l_0 = 0.07785m$. OF truss has a constant torque applied $\tau = 0.033 Nm$. Initial position is obtained from $\beta_0 = -0.062 rad$. Gravity is not considered. Using a classical position-based Newmark approach, the system fails with $\alpha = 0.5$ and $\beta = 0.25$. It has only been possible to obtain a solution with $\alpha = 0.5001$ and $\beta = 0.2507$. However, with the new approach, using $\alpha = 0,5$ and $\beta = 1/12$ (Fox

**Table 1** Coordinates for the fixed joints

| Joint | X | Y |
| --- | --- | --- |
| A | − 0.06934 | − 0.00227 |
| B | − 0.03635 | 0.03273 |
| C | 0.014 | 0.072 |
| O | 0 | 0 |

Goodwin), with a time increment of $2e-6s$ one reaches a good result, as presented in Fig. 24.

Although the aim of this paper is focused on stability rather than on precision, which obviously depends on the integrator used, it must be pointed out that the error obtained with Fox-Goodwin, with a $\Delta t = 2e-6s$ is under $2.28e-6$. The use of a quite good trapezoidal rule integrator (MBSLab), with a timestep of $1e-6s$ yields an error under $8.2e-6$. The use of Runge Kutta Fehlberg (an explicit method considered to be around 4$^{th}$–5th order, implemented in OpenSim), with a variable stepsize from $1e-5s$ to $1e-4s$ yields an error of $3.75e-7$. This seems to be reasonable taking into account that Fox-Goodwin is considered as a 4th order method and the trapezoidal rule is considered as 2nd order.

If one plots the maximum natural frequency of the system along the integration, Fig. 25 is obtained:

With this in mind, one could obtain a value for the stepsize that should not lead to instability. If one considers $\omega_{max} = 4503.477$, one can obtain from

Table 2 Mechanical properties of trusses

| Truss | Mass (kg) | Length (m) | $I_z$ (kg·m$^2$) | $X_G$ | $Y_G$ |
|---|---|---|---|---|---|
| OF | 0.004325 | 0.007 | $2.194 \cdot 10^{-6}$ | 0.00092 | 0 |
| EF | 0.00365 | 0.028 | $4.41 \cdot 10^{-7}$ | 0.0165 | 0 |
| HE | 0.00706 | 0.02 | $5.667 \cdot 10^{-7}$ | 0.00579 | 0 |
| GE | 0.00706 | 0.02 | $5.667 \cdot 10^{-7}$ | 0.00579 | 0 |
| AG | 0.05498 | 0.04 | $1.169 \cdot 10^{-5}$ | 0.02308 | 0.00916 |
| AH | 0.02373 | 0.04 | $1.912 \cdot 10^{-5}$ | 0.01228 | −0.00449 |

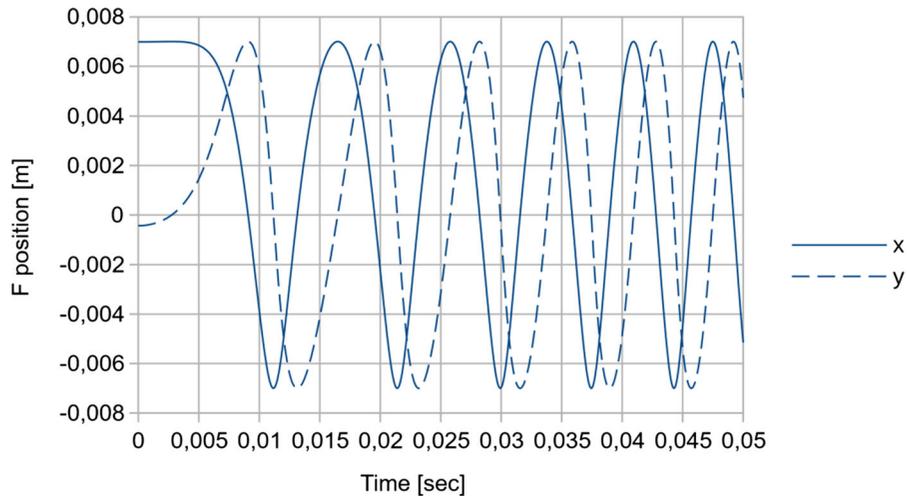

Fig. 24 Andrews Squeezer Mechanism Results using FG. Coordinates of Joint F

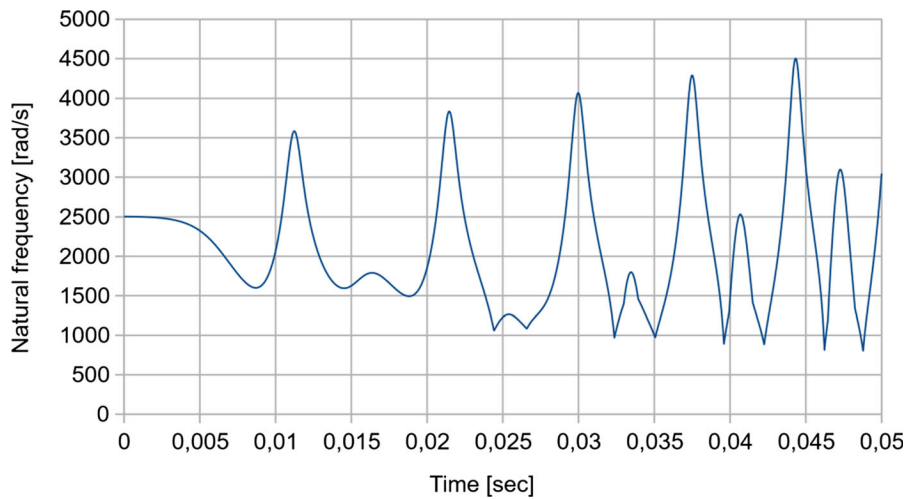

Fig. 25 Andrews Squeezer Mechanism. Evolution of maximum natural frequency

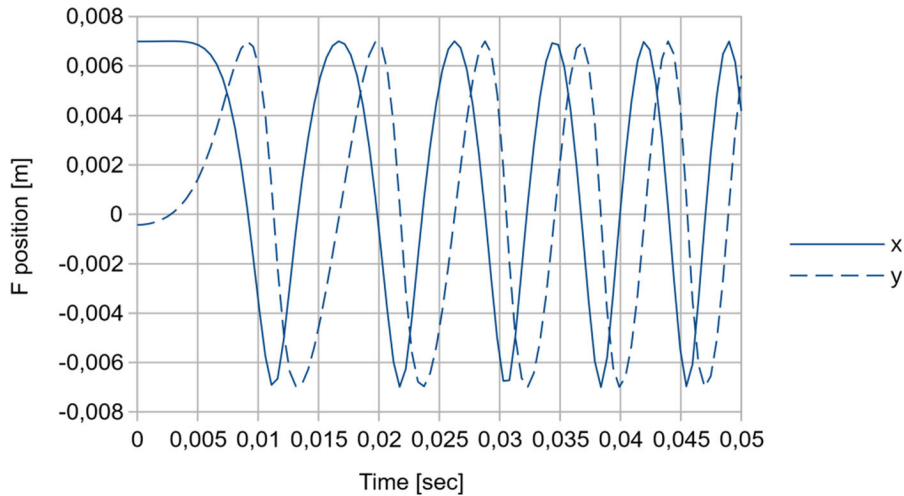

**Fig. 26** Andrews Squeezer Mechanism Results using FG. Coordinates of Joint F. Stepsize: $5e-4s$

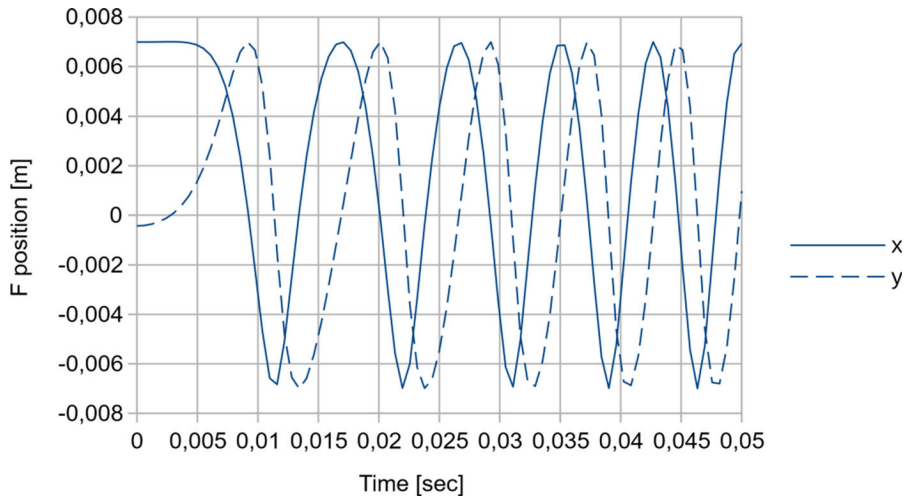

**Fig. 27** Andrews Squeezer Mechanism Results using FG. Coordinates of Joint F. Stepsize: $6e-4s$

Eq. (26) that using $\Delta_t \leq 5.44e - 4s$ verifies the linear condition for stability along all the integration time. The results obtained with $\Delta t = 5e - 4s$ are depicted in Fig. 26.

Evidently the error increases, but the system keeps being stable. A result that might surprise is that using $\Delta t = 6e - 4s$ (larger than the computed stability limit) the integration is also stable (see Fig. 27).

This is not actually completely unexpected. One must consider the fact that in nonlinear conditions natural frequencies vary over time and the computed stability limit is a worst-case scenario. The error might grow in a single iteration and in the following return to stability. It can even happen for some iterations and provoke no major stability issues. In fact, the system shows stability even with values of the stepsize around $1e - 3s$. In any case, it is obviously not a good practice to use a stepsize larger than the one that can be obtained from a stability analysis for the current integration step.

## 7 Some preliminary results regarding efficiency

In order to give an idea about the methods efficiency, a quite preliminary implementation in C has been carried out, which currently supports R and P joints. Although the implementation is quite preliminary, it has been developed taking into account performance and, thus, results should give an initial idea of the performance that one can expect from the method. The implementation is currently based in dense matrices.

A common problem to compare efficiency is to find benchmark results. Thanks to the IFTOMM Multibody Benchmark, some results are publicly available. Unfortunately no stiff problem is included there. Here we consider the term stiff in the sense of numerical integration and, therefore, related to the stability issues that the problem might lead to. This does not mean that the problem has no elastic elements. It means that the timestep required to correctly solve those problems is not large enough to lead to stability issues. Even the stiff flyball example is not stiff in the sense of numerical integration and can easily be solved by an explicit method. Taking into account that the stability does not affect efficiency, the performance can be measured with non-stiff problems. The studied examples are the double four bar and the rectangular Bricard mechanism. In both of them the benchmark limits the mechanical energy drift.

The first example is the double four bar mechanism. The maximum mechanical energy drift is 0.1 J. The published results are shown below (Table 3):

From the table above, one could state that a good result is of about 1 cent of a second, but one has to take those results carefully. Most of the results presented are provided using minimal, and relative coordinates, which heavily reduce the computational cost, but these methods are not general, and, thus, should not be compared with general methods as the one here presented. Also there are some 2D implementations that see their efficiency increased due to the reduction in the amount of variables. Thus, the results to be considered should be those presented by Tagliapietra (using OpenSim), Masarati (MBDyn), Urkullu (DIMCD, [46]) and Gonzalez (MBSLab). The reported times are, respectively, 0.456, 0.325, 0.0325 and 0.145. It is important to mention that Tagliapietra gets quite a good result in terms of precision, which hinders the real performance of the algorithm, which should be better. The method explained in this document takes 0.4 s (average of 10 simulations) on a computer with an Intel i7 7700HQ processor fixed at 2.8 GHz. Taking into account that this is an implicit, unconditionally stable method, it is a quite good result.

The second example to which this method has been applied is the Bricard mechanism. This system has redundant constraints. With this mechanism the aim is to keep the mechanical energy drift below 0.001 J. The results are shown below (Table 4):

Again, the results presented by Tagliapietra (using OpenSim) are misleading because of the use of a higher precision, which increases the cost. In any case, the times are 0.258 and 0.125. The here presented method takes 0.357 s in the same conditions exposed in the previous example. Again, a reasonable result.

**Table 3** Results of the double four bar mechanism

| Method | Author | Coord | Step (s) | e (J) | CPU (s) | Processor | GeekBench score |
|---|---|---|---|---|---|---|---|
| Index-3 ALF | **J. Cuadrado** | **Natural** | 1e-2 | **0.0917** | **0.6** | **C2duo E6550** | **1400** |
| Index-3 ALF | **A. Luaces** | **Relatives** | 1e-2 | **0.0137** | **1.7** | **C2duo E6550** | **1400** |
| MbsLab | **F. González** | **Natural 3D** | 1e-2 | **0.0917** | **0.145** | **C2duo E8400** | **1900** |
| OpenMBS | **R. Pastorino** | **Natural 2D** | 1e-2 | **0.0877** | **0.1278** | **BeagleBoneB** | **NA** |
| Index-3 ALF | **R. Pastorino** | **Natural 2D** | 1e-2 | **0.0877** | **0.0045** | **I7 3740QM** | **3800** |
| Non-recursive NEF | **M. Burkhardt** | **Mínimal** | **Variable** | **0.0015** | **0.0568** | **I7 3770** | **4000** |
| MBDyn | **P. Masarati** | **3D** | 8e-3 | **0.09** | **0.325** | **I7 2620 M** | **3000** |
| Non-recursive NEF | **M. Burkhardt** | **Mínimal** | 1e-2 | **0.0002** | **0.0191** | **SnapDrag S800** | **NA** |
| Non-recursive NEF | **M. Burkhardt** | **Mínimas** | 1e-2 | **0.0002** | **0.0023** | **I7 960** | **2700** |
| SolidWorks | **C. Chaojie** | **Natural 2D** | 1e-2 | **0.0001** | **48.77** | **C2Duo E8400** | **1900** |
| OpenSim | **L. Tagliapietra** | **Natural** | 1e-3 < > 1e-2 | **3.2e-7** | **0.455** | **I5 4570** | **4000** |
| Biolim | **F. Mouzo** | **Relatives 3D** | 1e-2 | **0.029** | **0.0226** | **i7 6700 K** | **5400** |
| DIMCD | **G. Urkullu** | **Cartesian 3D** | 1e-2 | **0.0885** | **0.0325** | **Ryzen 5 2600** | **4700** |

Table 4  Results of the Bricard mechanism

| Method | Author | Coords | Step (s) | e (J) | CPU (s) | Processor | GeekBench score |
|---|---|---|---|---|---|---|---|
| Index-3 ALF | F. González | Natural 3D | 5e-3 | 8e-4 | 0.125 | I7 3770 4000 | 4000 |
| OpenSim 3.2 | L. Tagliapietra | Natural | 1e-3 < > 1e-2 | 9.6e-7 | 0.258 | I5-4570 4100 | 4100 |
| DIMCD | G. Urkullu | Cartesian 3D | 1.5e-2 | 8.1e-4 | 0.0777 | Ryzen 5 2600 | 4700 |

## 8 Conclusions and future work

The integration of DAEs along with nonlinear constraints leads usually to unqualified DAEs. This translates into unpredictable stability behavior of the algorithms. This can be usually avoided by switching to minimal coordinates, but this arises the need of a particular mathematical development for each case, which can be an obstacle for general problems. The usual approach to solve this issue in general algorithms related to implicit methods is the use of heavily damped integrators such as Newmark with high values of the parameters or HHT, but this heavily impacts the convergence of the methods introducing a considerable amount of numerical damping. Here an alternative approach has been presented to perform the integration. In this approach, one integrates the DAE in a set of coordinates representing the tangent space of the manifold. To do this, a Taylor expansion of the constraints is iteratively performed in the current timestep which translates the function to a minimal coordinate set. This allows one to integrate in the tangent space. To obtain the initial values of the function in these minimal coordinates, these are projected in the tangent space. One advantage of this approach is that constraints are applied for each step simultaneously to the function to integrate and all the required derivatives. Furthermore, the new approach has not the unpredictable stability behavior of the classical approach, thus allowing one to use high order implicit methods such as Fox-Goodwin, or linear acceleration with predictable stability. Several tests have been performed with this algorithm, applying a Newmark approach, with Fox-Goodwin, linear acceleration and classical Newmark configurations. These confirm the predictable stability behavior and also yield to proper results. Currently no efficiency experimental tests have been performed. Due to the fact that, compared to the classical approach the new algorithm only includes operations of, at most, the same order ($N^3$) than the other methods, one can estimate that the new approach should not considerably increase computational cost, although it might not allow one to use certain sparse methods without modification. Additionally, a simple study for the case of structural integrators has been presented which explains the behavior of stability in DAEs when compared to unconstrained ODEs.

Future studies should be focused in performance. With this in mind a general Newmark approach should be developed in a high-performance compiled language. Additionally a deeper mathematical analysis of the method would be necessary. A proper analysis on the possible impact of the projection made for the values in the previous step in the convergence should also be performed.

**Acknowledgements**  The authors would like to thank to the Basque Government for its funding (ref. IT 947-16) and to the Spanish Ministry of Economy and Competitiveness for the grant through the project DPI2016-80372-R (AEI/FEDER, UE). We acknowledge Javier Cuadrado and his group for the help they have always given to us. We also thank Rafael Avilés for his mentorship.

**Funding**  Open Access funding provided thanks to the CRUE-CSIC agreement with Springer Nature.

**Declarations**

**Conflict of interest**  The authors declare that they have no conflict of interest.